\documentclass{article}
\usepackage{graphicx} 
\usepackage{amssymb,latexsym,amsfonts,amsthm,a4wide}              
\usepackage{mathrsfs}
\usepackage{mathtools}
\usepackage{bm}%
\usepackage{subfigure}
\usepackage{float}
\usepackage{cases}
\usepackage{textcomp} %
\usepackage{hyperref}
\usepackage{booktabs}%
\usepackage{color,xcolor}

\newtheorem{theorem}{Theorem}[section]
\newtheorem{lemma}{Lemma}[section]
\newtheorem{remark}{Remark}

\numberwithin{equation}{section}

\newcommand{\normmm}[1]{{\left\vert\kern-0.25ex\left\vert\kern-0.25ex\left\vert #1   
   \right\vert\kern-0.25ex\right\vert\kern-0.25ex\right\vert}}                       

\title{An Efficient Iterative Decoupling Method for Thermo-Poroelasticity Based on a Four-Field Formulation}

\author{Cai Mingchao\thanks{Morgan State University}
       \and  Li Jingzhi\thanks{South University of Science and Technology }  \and Li Ziliang \footnotemark[2]
       \and
       Liu Qiang\thanks{Shenzhen University}\footnotemark[3]
        }
\date{}

\begin{document}

\maketitle

\begin{abstract}
This paper studies the thermo-poroelasticity model. By introducing an intermediate variable, we transform the original three-field model into a four-field model. Building upon this four-field model, we present both a coupled finite element method and a decoupled iterative finite element method. We prove the stability and optimal convergence of the coupled finite element method. Furthermore, we establish the convergence of the decoupled iterative method. 
This paper focuses primarily on analyzing the iterative decoupled algorithm. It demonstrates that the algorithm's convergence does not require any additional assumptions about physical parameters or stabilization parameters.
Numerical results are provided to demonstrate the effectiveness and theoretical validity of these new methods.

\noindent\textbf{Keywords} thermo-poroelasticity, decoupled algorithms, finite element method.
\\

\end{abstract}

\maketitle

\section{Introduction} 
The thermo-poroelasticity model provides a comprehensive framework for describing the interplay between mechanical deformation, fluid flow, and heat transfer in porous media. By integrating Biot's poroelasticity theory with the heat equation, this model captures the effects of temperature variations on stress and displacement, accounting for the compressibility and thermal expansion of both fluid and solid components. The foundational work in poroelasticity was pioneered by Terzaghi and later expanded by Biot \cite{terzaghi1943theoretical, biot1941general, biot1955theory}, who developed a coupled model to describe solid-fluid and fluid-solid interactions. Extending these concepts to non-isothermal scenarios, the thermo-poroelasticity model comprises a force balance equation, a mass conservation equation, and a heat conduction equation, effectively representing the complex interactions among these physical processes. 
This model's broad applicability across various engineering and environmental disciplines has driven significant research interest. Applications span various areas such as soil mechanics, pharmaceutical tablet analysis, nuclear waste management \cite{yow2002coupled}, geothermal energy systems \cite{watanabe2011numerical}, biomedical engineering, and carbon dioxide sequestration \cite{rutqvist2002modeling}.

In recent years, there has been significant progress in computational research on the thermo-poroelasticity model. Studies such as \cite{brun2019well} have demonstrated the existence and uniqueness of solutions for nonlinear thermo-poroelasticity, along with deriving a priori energy estimates and regularity properties. To overcome computational challenges, the classical three-field formulation has been expanded to a four-field model by introducing the divergence of displacement as an auxiliary variable \cite{chen2022multiphysics}, which has also been extended to nonlinear thermo-poroelasticity, incorporating convective transport \cite{Ge2023analysis}. Various numerical methods have been proposed for these models, including Galerkin methods and combined Galerkin-mixed finite element methods for nonlinear formulations \cite{zhang2024coupling, zhang2022galerkin}. Additional approaches include a hybrid finite element method utilizing mixed finite elements for pressure, characteristic finite elements for temperature, and Galerkin finite elements for displacement \cite{Zhang2023mfe}, as well as discontinuous Galerkin methods and related algorithms for fully coupled nonlinear problems \cite{bonetti2024robust, antonietti2023discontinuous}.

Among recent numerical algorithms, iteratively decoupled algorithms have been extensively explored. For instance, unconditionally stable sequential methods for all-ways coupled thermo-poromechanics were proposed in \cite{Kim2018unconditionally}. Splitting schemes for decoupling poroelasticity and thermoelasticity \cite{Kolesov2014splitting}, and monolithic, partially decoupled, and fully decoupled schemes based on a five-field formulation incorporating heat and Darcy fluxes \cite{brun2020monolithic} have also been investigated. Furthermore, a decoupled iterative solution approach combined with reduced-order modeling has been developed to improve computational efficiency \cite{Ballarin2024projection}. While convergence proofs for coupled and decoupled iterative schemes exist \cite{Kolesov2014splitting, brun2020monolithic}, analysis of the convergence order of these methods remains limited.

This work focuses on developing an unconditionally stable and optimally convergent algorithm rooted in a novel four-field formulation of the thermo-poroelasticity model. To mitigate the issue of Poisson locking, we introduce an auxiliary variable \(\xi = -\lambda \nabla \cdot \bm{u} + \alpha p + \beta T\), transforming the original three-field model into a robust four-field framework. Building upon this reformulation, we design and rigorously analyze an iteratively decoupled method tailored for linear thermo-poroelastic problems. Theoretical analysis establishes the method's unconditional stability and optimal convergence. Unlike earlier approaches, such as time-extrapolation-based decoupling algorithms for Biot’s model \cite{ju2020parameter}, which often face stability limitations and reduced accuracy due to explicit extrapolation, our method overcomes these drawbacks. By iteratively solving each subproblem while incorporating time extrapolation, we ensure not only stability but also maintain optimal convergence rates. This approach provides a more reliable and accurate computational framework for addressing coupled thermo-poroelastic processes.

The structure of the paper is as follows:  
In Section \ref{sec: PDE model}, we provide an overview of the thermo-poroelastic problem. We introduce a four-field model and give it a priori estimates. 
In Section \ref{sec: algorithms}, we present two algorithms: the coupled algorithm and the decoupled algorithm.
In Section \ref{sec: convergence analysis}, we analyze the convergence of the decoupled algorithm. Furthermore, we derive the convergence rate for the discrete coupled problem relative to the continuous coupled problem.
Finally, Section \ref{sec: experiments} showcases the numerical results. 

\section{Mathematical model and its reformulation}\label{sec: PDE model}

We introduce some notations that will be utilized throughout the subsequent text.
Let $\Omega\in\mathbb R^d,~d=2,3$ be a bounded domain and $J=(0,\tau)$ with $\tau$ being the final time. 
For $1\le s < \infty$, let $L^s(\Omega)=\{v:\int_\Omega|v|^s<\infty \}$  represent the standard Banach space, with the associated norm $\Vert v\Vert_{L^s(\Omega)}=(\int_\Omega |v|^s)^{\frac1s}$.
We use $(\cdot,\cdot)$ to denote the $L^2$-inner product. 
For $s=\infty$, 
let $L^\infty(\Omega)=\{ v : \Omega \to \mathbb{R} : \text{ess sup}_{x \in \Omega} |v(x)| < \infty \} $ be the space of uniformly bounded measurable functions on $\Omega$,
with norm $\Vert v\Vert_{L^\infty(\Omega)}= \text{ess sup}_{x \in \Omega} |v|$.
Let $W_s^m(\Omega)=\{ v:D^\alpha v\in L^s(\Omega),0\le |\alpha|\le m,\|v\|\le \infty\}$ denote the standard Sobolev space with the associated norm $\|v\|_{W_s^k(\Omega)}=(\sum_{|\alpha|\le k}\|D^\alpha v\|^s_{L^s(\Omega)})^{1/s}$.
In particular, define $H^m(\Omega):= W_2^m(\Omega)$ as the special Sobolev space with $p=2$ and define
 $H_0^m(\Omega)=\{v\in H^1(\Omega):v|_{\partial \Omega}=0 \}$.
 For a Banach space $X$, we give $L^s(J;X)=\{v: (\int_J\|v\|_{X}^s)^{\frac1s}<\infty\}$.
For $v \in L^s(J;X)$, define the associated norm 
$\Vert v\Vert_{ L^s(J;X)}= (\int_J\|v\|_X^s)^{\frac1s}$.
For $v\in L^2(J,L^2(\Omega))$, it is a function in space $L^2(\Omega)$ square-integrable in time.  
Moreover, we define the space $L^\infty(J;L^2(\Omega))=\{v(t):\text{ess sup}_{t \in J}\|v(t)\|_{L^2(\Omega)}<\infty \}$.
We also use $\bm v_t$ and  $v_t$ to denote the partial derivative of the function $\bm v$ in $[H^1(\Omega)]^d$ and $ v$ in $H^1(\Omega)$ with respect to $t$, respectively.
For two positive quantities $x$ and $y$, we write $x \lesssim y$ to indicate that there exists a constant $C > 0$, depending only on the domain $\Omega$ and the final time $\tau$, such that $x \leq C y$. 

\subsection{A linear thermo-poroelasticity model}

We consider the following linear thermo-poroelasticity problem \cite{brun2018upscaling}, where the goal is to determine 
the displacement \(\bm{u}\), fluid pressure \(p\), and temperature \(T\). The governing equations are:  
\begin{equation}\label{TP_model}
\begin{aligned} 
    -\nabla \cdot (2\mu \bm{\varepsilon}(\bm{u}) + \lambda \nabla \cdot \bm{u} \bm{I}) + \alpha \nabla p + \beta \nabla T &= \bm{f}, \quad \text{in } \Omega \times J, \\
    \frac{\partial}{\partial t}(c_0 p - b_0 T + \alpha \nabla \cdot \bm{u}) - \nabla \cdot (\bm{K} \nabla p) &= g, \quad \text{in } \Omega \times J, \\
    \frac{\partial}{\partial t}(a_0 T - b_0 p + \beta \nabla \cdot \bm{u}) - \nabla \cdot (\bm{\Theta} \nabla T) &=  H_{s}, \quad \text{in } \Omega \times J,
\end{aligned}
\end{equation}

Proper boundary and initial conditions should be provided in order to ensure the
existence and uniqueness of the solution.  In this paper, we consider a mixed partial
Neumann and partial Dirichlet conditions for $\bm{u}$ and the Dirichlet conditions for $p$ and $T$. Assume $\partial \Omega =\Gamma_d+\Gamma_n$ with
$\Gamma_d \neq \partial\Omega $.  Here, $\Gamma_d$ and $\Gamma_n$ denote the Dirichlet boundary and Neumann boundary for $\bm{u}$, respectively. Without loss of generality, the boundary conditions are assumed to be homogeneous.
Specifically, the boundary conditions are given by
\begin{equation}\label{eq: B_C}
\begin{aligned}
    \bm{u} = 0, \quad& \text{on } \Gamma_d,\\
     ( 2\mu \bm{\varepsilon}(\bm{u}) + \lambda \nabla \cdot \bm{u} \bm{I} - \alpha  p\bm I - \beta  T \bm I)\bm n = 0 , \quad& \text{on } \Gamma_n   ,    \\
    p=0, \quad& \text{on } \partial\Omega,  \\
    T=0, \quad& \text{on } \partial\Omega.
\end{aligned}
\end{equation}
Here, $\bm{n}$ is the unit outward normal to the boundary.
 The initial conditions are given by 
 \begin{equation}\label{eq: I_C}
\bm{u}(\cdot, 0) = \bm{u}^0, \quad p(\cdot, 0) = p^0, \quad T(\cdot, 0) = T^0.
\end{equation}

The operator $\bm\varepsilon$ is defined as $\bm\varepsilon(\bm u)=\frac{1}{2}(\nabla\bm u+\nabla\bm u^T)$ and
the matrix $\bm I$  is the identity tensor. The right-hand side functions $\bm f,g$, and $ H_{s}$ denote the body force, the mass source, and the heat source, respectively. $\alpha$ is the Biot-Willis constant and $\beta$ is the thermal stress coefficient. $a_0,b_0$, and $c_0$ represent the effective thermal capacity, the thermal dilation coefficient, and the specific storage coefficient, respectively. The matrix parameter $\bm K = (K_{ij})^d_{ij=1}$ 
denotes the permeability divided by fluid viscosity and matrix parameter  $\bm\Theta = (\Theta_{ij})^d_{ij=1}$
denotes the effective thermal conductivity. $\mu$ and $\lambda$ represent the Lam\'e parameters, which can be formulated in terms of Young's modulus $E$ and Poisson's ratio $\nu$:
\[
\lambda=  \frac{E\nu}{(1+\nu)(1-2\nu)},\qquad \mu=\frac{E}{2(1+\nu)}.
\]
Note that in this paper, the considered coefficients and variables in the model are without units, which is also known as the dimensionless form.

We now introduce the four-field formulation for the above model \eqref{TP_model}. More clearly, following the methodology for handling Biot's problem in \cite{lee2017parameter, oyarzua2016locking}, we define an intermediate auxiliary variable to represent the volumetric contribution to the total stress:  
$$  
\xi = -\lambda \nabla \cdot \bm{u} + \alpha p + \beta T,
$$  
which is commonly referred to as the pseudo-total pressure \cite{antonietti2023discontinuous}. Substituting this variable into equation \eqref{TP_model} and using parameter transformations above, the four-field thermo-poroelasticity problem can be expressed as:

\begin{equation}\label{TP_Model_four}
\begin{aligned} 
       -2\mu \nabla\cdot \bm\varepsilon(\bm u) +\nabla\xi 
                      &=\bm f, \\
      -\lambda \nabla\cdot\bm u-\xi+\alpha p+\beta T
                      &=0, \\
-\frac\alpha\lambda\frac\partial{\partial t}\xi+(c_0+\frac{\alpha^2}\lambda)\frac\partial{\partial t}p
                 -\nabla\cdot(\bm K\nabla p)+(\frac{\alpha\beta}\lambda-b_0)\frac\partial{\partial t}T
                      &=g, \\
-\frac\beta\lambda\frac\partial{\partial t}\xi +(\frac{\alpha\beta}\lambda-b_0)\frac\partial{\partial t}p
+(a_0+\frac{\beta^2}\lambda)\frac\partial{\partial t}T
                  -\nabla\cdot(\bm \Theta\nabla T)
                      &= H_{s}. 
\end{aligned}
\end{equation}
We still use the initial condition (\ref{eq: I_C}) and the Dirichlet boundary conditions (\ref{eq: B_C}) for \eqref{TP_Model_four}. In practical situations, nonhomogeneous Dirichlet and Neumann boundary conditions are commonly encountered. The analysis performed for homogeneous boundary conditions can be straightforwardly extended to accommodate these cases.

As outlined in \cite{brun2020monolithic}, we make the following \textbf{Assumptions} for the model parameters throughout this paper (similar assumptions are also discussed in \cite{chen2022multiphysics, zhang2024coupling}):
 
(A1) $\bm K$ and $\bm\Theta$ are constant in time, symmetric, definite, and positive. There exist $k_m>0$
and $k_M>0$ such that
$$ k_m|\bm v|^2 \le \bm v^T\bm K(x)\bm v,\quad |\bm K(x)\bm v|\le k_M|\bm v|,\quad \forall \bm v\in\mathbb R^d \backslash \{\bm 0\}, $$
where $|\cdot|$ denotes the $l^2$-norm of a vector in $\mathbb R^d$.
There exist  $\theta_m>0$ and $\theta_M>0$ such that
$$ \theta_m|\bm v|^2 \le \bm v^T\bm\Theta(x)\bm v,\quad |\bm\Theta(x)\bm v|\le\theta_M|\bm v|,\quad \forall \bm v\in\mathbb R^d 
\backslash \{\bm 0\}. $$

(A2) The constants \(\lambda\), \(\mu\), \(\alpha\) and \(\beta\) are strictly positive constants.

(A3) The constants $a_0,b_0 $ and $c_0$ are such that $a_0,~ c_0 \ge b_0 \ge 0$ .

In addition to \textbf{Assumptions} (A1)-(A3), the following Assumptions will also be maintained throughout the remainder of this paper. 

We assume the initial conditions satisfy \(\bm{u}^0 \in [H^1(\Omega)]^d\) and \(p^0, T^0 \in L^2(\Omega)\), while the source terms satisfy \(\bm{f} \in [L^2(\Omega)]^d\), \(g \in L^2(\Omega)\), and \( H_{s} \in L^2(\Omega)\). For simplicity, we further assume that \(\bm{f}\), \(g\), and \( H_{s}\) are time-independent. Define the following spaces
$$
\bm V=\bm H^1_{0,\Gamma_d}(\Omega),~Q=L^2(\Omega),~W=H_0^1(\Omega). 
$$
and their dual spaces are denoted as  $\bm V'$,~$Q'$ and $W'$. 
It is a well-established result that the following inf-sup condition holds: there exists a constant $\zeta_0>0$ depending only on $\Omega$ and $\Gamma_d$ such that 
\begin{equation}\label{infsup V Q}
   \zeta_0
         \le \inf_{\xi\in L^2(\Omega)} \sup_{\bm v\in\bm V}\frac{|(\text{div }\bm v,\xi)|}{\|\bm v\|_{H^1(\Omega)}\| \xi\|_{L^2(\Omega)}}.
\end{equation}
Its proof can be found in Quarteroni's book \cite{Quarteroni1999} or Cai's thesis \cite{cai2008modeling} (Lemma 2.2). Given $\tau>0$, a 4-tuple $(\bm u,\xi,p,T)$ with
$$
\bm u\in L^\infty(J;\bm V),~\xi\in L^\infty(J; Q),~p,~T\in L^\infty(J; W), 
$$
$$
\nabla u_t, ~\xi_t, ~p_t, ~T_t\in L^2(J;  W').
$$
is called a weak solution of problem \eqref{TP_Model_four}, if for almost every $t\in (0, \tau]$, and $\forall (\bm v,~\phi,~q,~S)\in\bm V\times Q\times W\times W$, there holds
\begin{equation}\label{TP_Model_var}
\begin{aligned} 
2\mu(\bm\varepsilon(\bm u), \bm\varepsilon(\bm v)) - (\nabla\cdot\bm v,\xi)&=(\bm f,\bm v)       ,~\forall \bm v\in\bm V,  \\
  -(\nabla\cdot\bm u,\phi)-\frac{1}{\lambda}(\xi,\phi)+\frac{\alpha}{\lambda}(p,\phi)+\frac{\beta}{\lambda}(T,\phi)&=0 ,~\forall\phi\in Q,\\
 -\frac\alpha\lambda( \xi_t ,q)+(c_0+\frac{\alpha^2}{\lambda})( p_t ,q)+\bm (\bm K\nabla p,\nabla q)+ (\frac{\alpha\beta}{\lambda}-b_0)( T_t ,q)
                &=(g,q),~\forall q\in W,\\
 -\frac\beta\lambda( \xi_t ,S)+(\frac{\alpha\beta}{\lambda}-b_0)( p_t ,S)+(a_0+\frac{\beta^2}{\lambda})( T_t ,S)   
  +(\bm\Theta\nabla T,\nabla S)&=( H_{s},S), ~\forall S\in W.\\
\end{aligned}
\end{equation}  
Here, these equations can be derived by multiplying the differential equations with appropriate test functions. 
By the definition of  $\xi = -\lambda \nabla \cdot \bm{u} + \alpha p + \beta T$, we have 
\begin{equation}\label{eq: initial xi0}
 \xi(\cdot,0) = \xi^0 = -\lambda \nabla \cdot \bm{u}^0 + \alpha p^0 + \beta T^0. 
\end{equation}
Then, the initial condition \eqref{eq: I_C}, \eqref{eq: initial xi0} and boundary condition \eqref{eq: B_C} can still be applied to the
Problem \eqref{TP_Model_var}.

\subsection{A priori estimates of the weak solution}

In this subsection, we provide an energy estimate and a priori estimates. Using these estimates, we establish the existence and uniqueness theorem in Theorem \ref{Th: exist and uniq}. Firstly, we present Korn's inequality(\cite{brenner2004korn}), which elucidates the relationship between \(\|\bm u\|_{H^1(\Omega)}\) and \(\|\bm\varepsilon(\bm u)\|_{L^2(\Omega)}\). 
\begin{lemma}\label{Korn's inequa}
There exists a constant $C$ such that Korn’s inequality 
\begin{equation*}
               C \Vert  \bm  u\|_{H^1(\Omega)}^2   \le  \Vert\bm\varepsilon(\bm u)\|_{L^2(\Omega)}^2,
               \quad\forall \bm u\in [H^1_{0,\Gamma_d}(\Omega)]^d,~d=2,3
\end{equation*}
holds true.
\end{lemma}
Then, we give an energy estimate.
 \begin{lemma}\label{le: energy law}
 Assume that the {\bf Assumptions} of (A1)-(A3) hold. 
 Every weak solution $(\bm u,\xi,p,T)$ of problem \eqref{TP_Model_var} satisfies the following energy law:
\begin{equation}\label{EtE0 energy law}
\begin{aligned}
E(t)+\int_0^t \bm K(\nabla p,\nabla p)ds+\int_0^t\bm\Theta(\nabla T,\nabla T)ds
-\int_0^t(g,p)ds-\int_0^t( H_{s},T)ds
=E(0),
\end{aligned}
\end{equation}  
for all $t\in J$, where 
 \begin{equation}\label{eq: E(t)}
 \begin{aligned}
      E(t):=&\mu\Vert\bm\varepsilon(\bm u)\Vert_{L^2(\Omega)}^2+\frac1{2\lambda}\Vert\alpha p+\beta T-\xi\Vert_{L^2(\Omega)}^2
  +\frac{c_0-b_0}{2}\Vert p\Vert_{L^2(\Omega)}^2\\
  &+\frac{a_0-b_0}{2}\Vert T\Vert_{L^2(\Omega)}^2+\frac{b_0}{2}\Vert p-T\Vert_{L^2(\Omega)}^2-(\bm f,\bm u).
 \end{aligned}
\end{equation}
Furthermore, 
 \begin{equation}\label{eq: bound xi}
 \begin{aligned}
  \|\xi(t)\|_{L^2(\Omega)}\le C\big(  2\mu\|\bm \varepsilon(\bm u)\|_{L^2(\Omega)}   +\|\bm f\|_{L^2(\Omega)}   \big),
 \end{aligned}
\end{equation} 
where $C$ is a constant depending only on $\Omega$ and $\Gamma_d$.
\end{lemma}
\begin{proof}
 Differentiating the second equation of system \eqref{TP_Model_var} with respect to $t$, taking $\bm v= \bm u_t,\phi=-\xi,q=p$ and $S=T$ in \eqref{TP_Model_var}, we have
\begin{equation*}\label{TP_Model_var_eq2dt}
\begin{aligned} 
2\mu(\bm\varepsilon(\bm u),\bm\varepsilon( \bm u_t)) - (\nabla\cdot \bm u_t,\xi)&=(\bm f, \bm u_t)       ,   \\
  (\nabla\cdot \bm u_t,\xi)+\frac{1}{\lambda}( \xi_t,\xi)-\frac{\alpha}{\lambda}( p_t,\xi)-\frac{\beta}{\lambda}( T_t,\xi)&=0 ,\\
 -\frac\alpha\lambda( \xi_t ,p)+(c_0+\frac{\alpha^2}{\lambda})( p_t ,p)+  (\bm K\nabla p,\nabla p)+ (\frac{\alpha\beta}{\lambda}-b_0)( T_t ,p)
                &=(g,p),\\
 -\frac\beta\lambda( \xi_t ,T)+(\frac{\alpha\beta}{\lambda}-b_0)( p_t ,T)+(a_0+\frac{\beta^2}{\lambda})( T_t ,T)   
  +(\bm\Theta\nabla T,\nabla T)&=( H_{s},T).\\
\end{aligned}
\end{equation*}
Adding the above four equations together, we have
\begin{equation}\label{summation utxtptTt}
\begin{aligned} 
&2\mu(\bm\varepsilon(\bm u),\bm\varepsilon( \bm u_t))+\frac{1}{\lambda}( \xi_t,\xi)-\frac{\alpha}{\lambda}( p_t,\xi)-\frac{\beta}{\lambda}( T_t,\xi)\\
& -\frac\alpha\lambda( \xi_t ,p)+(c_0+\frac{\alpha^2}{\lambda})( p_t ,p)+(\bm K \nabla p,\nabla p)+ (\frac{\alpha\beta}{\lambda}-b_0)( T_t ,p)\\
& -\frac\beta\lambda( \xi_t ,T)+(\frac{\alpha\beta}{\lambda}-b_0)( p_t ,T)+(a_0+\frac{\beta^2}{\lambda})( T_t ,T)   
  +(\bm\Theta\nabla T,\nabla T)\\
& -(\bm f, \bm u_t) -(g,p) -( H_{s},T)=0,\\
\end{aligned}
\end{equation}

From the expression of \eqref{eq: E(t)}, we note that
\begin{equation} \label{dt the Et}
\begin{aligned} 
  &\frac{d}{dt}\Big( \frac1{2\lambda}\Vert\alpha p+\beta T-\xi\Vert_{L^2(\Omega)}^2
  +\frac{c_0-b_0}{2}\Vert p\Vert_{L^2(\Omega)}^2
  +\frac{a_0-b_0}{2}\Vert T\Vert_{L^2(\Omega)}^2+\frac{b_0}{2}\Vert p-T\Vert_{L^2(\Omega)}^2
 \Big)\\
 =&\frac1\lambda(\alpha p+\beta T-\xi,\alpha p_t+\beta T_t- \xi_t)
 +(c_0-b_0)(p, p_t)+(a_0-b_0)(T, T_t)\\
 &+b_0(p-T, p_t- T_t)
 \\
 =&\frac{1}{\lambda}( \xi_t,\xi)-\frac{\alpha}{\lambda}( p_t,\xi)-\frac{\beta}{\lambda}( T_t,\xi) 
  -\frac\alpha\lambda( \xi_t ,p)
  +(c_0+\frac{\alpha^2}{\lambda})( p_t ,p) \\
&  + (\frac{\alpha\beta}{\lambda}-b_0)( T_t ,p) 
  -\frac\beta\lambda( \xi_t ,T)
  +(\frac{\alpha\beta}{\lambda}-b_0)( p_t ,T)\\
&  +(a_0+\frac{\beta^2}{\lambda})( T_t ,T).   
  \\ 
\end{aligned}
\end{equation} 
Integrating \eqref{summation utxtptTt} from $0$ to $t$ and using \eqref{dt the Et}, we derive \eqref{EtE0 energy law}.
The bound for $\xi$ follows from the inf-sup condition between $\bm V$ and $Q$ in \eqref{infsup V Q} and the Korn's inequality. We have the following inequality holds:
\begin{equation*}\label{eq:bound xi proof}
\begin{aligned} 
\zeta_0\| \xi\|_{L^2(\Omega)}
         \le&\sup_{\bm v\in\bm V}\frac{|(\text{div }\bm v,\xi)|}{\|\bm v\|_{H^1(\Omega)}}\\
  \le&\sup_{\bm v\in\bm V}\frac{ | 2\mu(\bm\varepsilon(\bm u),\bm\varepsilon(\bm v))| +|(\bm f,\bm  v)|}
           {\|\bm v\|_{H^1(\Omega)}}\\
  \le& C_K\big(  2\mu\|\bm \varepsilon(\bm u)\|_{L^2(\Omega)}   +\|\bm f\|_{L^2(\Omega)}   \big),      
\end{aligned}
\end{equation*} 
where constant $\zeta_0$ is from the inf-sup condition  and $C_k$ is from the Korn’s
inequality. This deduces \eqref{eq: bound xi}. 
\end{proof}

Then, the energy law \ref{le: energy law}  implies the following priori estimates.
\begin{theorem}\label{Th: boundness1}
Assume that the {\bf Assumptions} of (A1)-(A3) hold. 
There exists a positive constant
 $C=C(\Vert\bm\varepsilon(\bm u^0)\Vert_{L^2(\Omega)}$, $\Vert\xi^0\Vert_{L^2(\Omega)}$, $\Vert p^0\Vert_{L^2(\Omega)}$, $\Vert T^0\Vert_{L^2(\Omega)}$, $\Vert \bm f\Vert_{L^2(\Omega)}$, $
\Vert g \Vert_{L^2(\Omega)}$, $\Vert  H_{s}\Vert_{L^2(\Omega)})$
 such that
 \begin{equation}\label{estimate C}
\begin{aligned} 
 &\sqrt{2\mu}\Vert\bm\varepsilon(\bm u)\Vert_{L^\infty(J;L^2(\Omega))}
 +\sqrt{\frac{1}{\lambda}}\Vert\alpha p+\beta T-\xi\Vert_{L^\infty(J;L^2(\Omega))}+\sqrt{c_0-b_0}\Vert p\Vert_{L^\infty(J;L^2(\Omega))}\\
 &+\sqrt{a_0-b_0}\Vert T\Vert_{L^\infty(J;L^2(\Omega))} +\sqrt{b_0}\Vert p-T\Vert_{L^\infty(J;L^2(\Omega))}
 +\sqrt{k_m}\Vert \nabla p\Vert_{L^2(J;L^2(\Omega))}\\
 &+\sqrt{\theta_m}\Vert \nabla T\Vert_{L^2(J;L^2(\Omega))}
 \leq C.
 \end{aligned}
\end{equation} 
\end{theorem}
\begin{proof}
    From \eqref{summation utxtptTt},  by using  the \textbf{Assumption} (A1), the Cauchy-Schwarz inequality, and Young's inequality, we have
\begin{equation}\label{Et le fut+gp+ht}
\begin{aligned} 
& \mu\frac{d}{dt}\Vert\bm\varepsilon(\bm u)\Vert_{L^2(\Omega)}^2  
  +\frac1{2\lambda}\frac{d}{dt}\Vert\alpha p+\beta T-\xi\Vert_{L^2(\Omega)}^2
  +\frac{c_0-b_0}{2}\frac{d}{dt}\Vert p\Vert_{L^2(\Omega)}^2  
 +\frac{a_0-b_0}{2}\frac{d}{dt}\Vert T\Vert_{L^2(\Omega)}^2 \\
&+\frac{b_0}{2}\frac{d}{dt}\Vert p-T\Vert_{L^2(\Omega)}^2  
+ k_m\Vert\nabla p\Vert_{L^2(\Omega)}^2
 +\theta_m\Vert\nabla T\Vert_{L^2(\Omega)}^2\\
&\le   (\bm f, \bm u_t) + (g,p) + ( H_{s},T).\\
\end{aligned}
\end{equation}
Integrating \eqref{Et le fut+gp+ht} from $0$ to $t$, applying the Cauchy-Schwarz inequality, Young's inequality, we can find a constant $\mu^*$ such that
 \begin{equation*}    
\begin{aligned} 
& \mu \Vert\bm\varepsilon(\bm u)\Vert_{L^2(\Omega)}^2  
  +\frac1{2\lambda} \Vert\alpha p+\beta T-\xi\Vert_{L^2(\Omega)}^2
  +\frac{c_0-b_0}{2} \Vert p\Vert_{L^2(\Omega)}^2  
 +\frac{a_0-b_0}{2} \Vert T\Vert_{L^2(\Omega)}^2 \\
&+\frac{b_0}{2} \Vert p-T\Vert_{L^2(\Omega)}^2  
+ \int_0^t k_m\Vert\nabla p\Vert_{L^2(\Omega)}^2ds
 +\int_0^t \theta_m\Vert\nabla T\Vert_{L^2(\Omega)}^2ds\\
&\le \mu \Vert\bm\varepsilon(\bm u^0)\Vert_{L^2(\Omega)}^2  
  +\frac1{2\lambda} \Vert\alpha p^0+\beta T^0-\xi^0\Vert_{L^2(\Omega)}^2
  +\frac{c_0-b_0}{2} \Vert p^0\Vert_{L^2(\Omega)}^2  
 +\frac{a_0-b_0}{2} \Vert T^0\Vert_{L^2(\Omega)}^2 \\
&+\frac{b_0}{2} \Vert p^0-T^0\Vert_{L^2(\Omega)}^2   
 +  2\mu^* \Vert \bm f\Vert_{L^2(\Omega)}^2
  +\frac1{4\mu^*}   \Vert\bm u\Vert_{L^2(\Omega)}^2 
  +\frac1{4\mu^*}   \Vert\bm u^0\Vert_{L^2(\Omega)}^2 \\
& + \frac12\int_0^t(\Vert g\Vert_{L^2(\Omega)}^2+\Vert  H_{s}\Vert_{L^2(\Omega)}^2+\Vert p\Vert_{L^2(\Omega)}^2+\Vert T\Vert_{L^2(\Omega)}^2) .\\
\end{aligned}
\end{equation*}
Then, by the Poincar\'{e} inequality, Gronwall's inequality, and Korn's inequality, choosing $\mu^*=\frac1{2\mu C_K}$, we get
  \begin{equation*}    
\begin{aligned} 
& \mu \Vert\bm\varepsilon(\bm u)\Vert_{L^2(\Omega)}^2  
  +\frac1{2\lambda} \Vert\alpha p+\beta T-\xi\Vert_{L^2(\Omega)}^2
  +\frac{c_0-b_0}{2} \Vert p\Vert_{L^2(\Omega)}^2  
 +\frac{a_0-b_0}{2} \Vert T\Vert_{L^2(\Omega)}^2 \\
&+\frac{b_0}{2} \Vert p-T\Vert_{L^2(\Omega)}^2  
+ \int_0^t k_m\Vert\nabla p\Vert_{L^2(\Omega)}^2ds
 +\int_0^t \theta_m\Vert\nabla T\Vert_{L^2(\Omega)}^2ds\\
&\lesssim \frac{3\mu}{2} \mu \Vert\bm\varepsilon(\bm u^0)\Vert_{L^2(\Omega)}^2  
  +\frac1{2\lambda} \Vert\alpha p^0+\beta T^0-\xi^0\Vert_{L^2(\Omega)}^2
  +\frac{c_0-b_0}{2} \Vert p^0\Vert_{L^2(\Omega)}^2  
 +\frac{a_0-b_0}{2} \Vert T^0\Vert_{L^2(\Omega)}^2 \\
&+\frac{b_0}{2} \Vert p^0-T^0\Vert_{L^2(\Omega)}^2   
 + \frac1{\mu} \Vert \bm f\Vert_{L^2(\Omega)}^2
  +\frac\mu{2} \Vert\bm\varepsilon(\bm u)\Vert_{L^2(\Omega)}^2  
  + \frac12\int_0^t(\Vert g\Vert_{L^2(\Omega)}^2+\Vert  H_{s}\Vert_{L^2(\Omega)}^2) ,\\
\end{aligned}
\end{equation*}
 which implies that \eqref{estimate C} holds.
\end{proof}


 

 %

 \begin{lemma}\label{Th: boundness2}
Assume that the {\bf Assumptions} of (A1)-(A3) hold. Then, there exist positive constants 
 $C_1=C_1(\Vert\bm\varepsilon(\bm u^0)\Vert_{L^2(\Omega)},~\Vert\xi^0\Vert_{L^2(\Omega)},~\Vert p^0\Vert_{L^2(\Omega)},~\Vert T^0\Vert_{L^2(\Omega)}$, $\Vert \bm f\Vert_{L^2(\Omega)}$, $
\Vert g \Vert_{L^2(\Omega)},\Vert  H_{s}\Vert_{L^2(\Omega)}$, $\Vert\nabla p(0)\Vert_{L^2(\Omega)}$, $\Vert\nabla T(0)\Vert_{L^2(\Omega)})$
 and
  $C_2=C_2( C_1,\Vert\bm u^0\Vert_{H^2(\Omega)},~\Vert p^0\Vert_{H^2(\Omega)},~\Vert T^0\Vert_{H^2(\Omega)})$
 such that
  \begin{equation}\label{ut2+alpha beta C1}
\begin{aligned} 
 &\sqrt{2\mu}\Vert\bm\varepsilon(\bm u_t)\Vert_{L^2(J;L^2(\Omega))}
 +\sqrt{\frac{1}{\lambda}}\Vert\alpha p_t+\beta T_t-\xi_t\Vert_{L^2(J;L^2(\Omega))}+\sqrt{c_0-b_0}\Vert p_t\Vert_{L^2(J;L^2(\Omega))}\\
 &+\sqrt{a_0-b_0}\Vert T_t\Vert_{L^2(J;L^2(\Omega))} +\sqrt{b_0}\Vert p_t-T_t\Vert_{L^2(J;L^2(\Omega))}
 +\sqrt{k_m}\Vert \nabla p_t\Vert_{L^\infty(J;L^2(\Omega))}\\
 &\sqrt{\theta_m}\Vert \nabla T_t\Vert_{L^\infty(J;L^2(\Omega))}
 \leq C_1,
 \end{aligned}
\end{equation}
and
 \begin{equation}\label{ut2+alpha beta C2}
\begin{aligned} 
 &\sqrt{2\mu}\Vert\bm\varepsilon(\bm u_t)\Vert_{L^\infty(J;L^2(\Omega))}
 +\sqrt{\frac{1}{\lambda}}\Vert\alpha p_t+\beta T_t-\xi_t\Vert_{L^\infty(J;L^2(\Omega))}+\sqrt{c_0-b_0}\Vert p_t\Vert_{L^\infty(J;L^2(\Omega))}\\
 &+\sqrt{a_0-b_0}\Vert T_t\Vert_{L^\infty(J;L^2(\Omega))}+\sqrt{b_0}\Vert p_t-T_t\Vert_{L^\infty(J;L^2(\Omega))}
 +\sqrt{2k_m}\Vert \nabla p_t\Vert_{L^2(J;L^2(\Omega))}\\
 &\sqrt{2\theta_m}\Vert \nabla T_t\Vert_{L^2(J;L^2(\Omega))} \leq C_2.
 \end{aligned}
\end{equation}
\end{lemma}

 \begin{proof}
Differentiating the first and the second equation of system \eqref{TP_Model_var} with respect to $t$,
 taking $\bm v = \bm u_t ,~\phi = - \xi_t,~q = p_t$ , and $S = T_t$  respectively, and adding the resulting equations together we have
\begin{equation*}    
\begin{aligned} 
&2\mu(\bm\varepsilon( \bm u_t),\bm\varepsilon( \bm u_t))+\frac{1}{\lambda}( \xi_t, \xi_t)-\frac{\alpha}{\lambda}( p_t, \xi_t)-\frac{\beta}{\lambda}( T_t, \xi_t)\\
& -\frac\alpha\lambda( \xi_t , p_t)+(c_0+\frac{\alpha^2}{\lambda})( p_t , p_t)+(\bm K \nabla p,\nabla p_t)+ (\frac{\alpha\beta}{\lambda}-b_0)( T_t , p_t)\\
& -\frac\beta\lambda( \xi_t , T_t)+(\frac{\alpha\beta}{\lambda}-b_0)( p_t , T_t)+(a_0+\frac{\beta^2}{\lambda})( T_t , T_t)   
  +(\bm\Theta\nabla T,\nabla T_t)\\
&=   (g, p_t) + ( H_{s}, T_t).\\
\end{aligned}
\end{equation*}
 Based on some algebraic operations, we have
\begin{equation}\label{norm inequality km thetam C1}
\begin{aligned} 
& 2\mu\Vert\bm\varepsilon( \bm u_t)\Vert_{L^2(\Omega)}^2  
  +\frac1{ \lambda}\Vert\alpha  p_t+\beta  T_t- \xi_t\Vert_{L^2(\Omega)}^2
  +(c_0-b_0) \Vert p_t\Vert_{L^2(\Omega)}^2  
 +(a_0-b_0) \Vert T_t\Vert_{L^2(\Omega)}^2 \\
&+ b_0  \Vert  p_t- T_t\Vert_{L^2(\Omega)}^2  
+(\bm K \nabla p,\nabla p_t)
 +(\bm\Theta\nabla T,\nabla T_t)\\
&=    (g, p_t) + ( H_{s}, T_t).\\
\end{aligned}
\end{equation}
Integrating equation \eqref{norm inequality km thetam C1} from $0$ to $t$, by {\bf Assumption} (A1) we have
\begin{equation*}  
\begin{aligned} 
& \int_0^t\Big( 2\mu\Vert\bm\varepsilon( \bm u_t)\Vert_{L^2(\Omega)}^2  
  +\frac1{ \lambda}\Vert\alpha  p_t+\beta  T_t- \xi_t\Vert_{L^2(\Omega)}^2
  +(c_0-b_0) \Vert p_t\Vert_{L^2(\Omega)}^2  
 +(a_0-b_0) \Vert T_t\Vert_{L^2(\Omega)}^2 \\
&+ b_0  \Vert  p_t- T_t\Vert_{L^2(\Omega)}^2 \Big)ds 
+ k_m \Vert\nabla p\Vert_{L^2(\Omega)}^2
 +\theta_m \Vert\nabla T\Vert_{L^2(\Omega)}^2\\
&\le    (g,p(t)-p(0)) + ( H_{s},T(t)-T(0)) 
      + k_M \Vert\nabla p(0)\Vert_{L^2(\Omega)}^2
 +\theta_M \Vert\nabla T(0)\Vert_{L^2(\Omega)}^2,\\
\end{aligned}
\end{equation*}
which, when combined with equation \eqref{estimate C} in lemma \ref{le: energy law} and Gronwall's inequality, implies that \eqref{ut2+alpha beta C1} holds.

Differentiating the first, the third, and the fourth equations of \eqref{TP_Model_var} with respect to $t$,
 taking $\bm v = \bm u_{tt} ,\phi = - \xi_t, q = p_t$ , and $S = T_t$  respectively, and adding the resulting equations we have
\begin{equation*}    
\begin{aligned} 
&2\mu(\bm\varepsilon( \bm u_{tt}),\bm\varepsilon( \bm u_t))
 +\frac{1}{\lambda}(  \xi_{tt}, \phi_t)
 -\frac{\alpha}{\lambda}(  p_{tt}, \phi_t)
 -\frac{\beta}{\lambda}(  T_{tt}, \phi_t)\\
& -\frac\alpha\lambda(  \xi_{tt} , p_t)+(c_0+\frac{\alpha^2}{\lambda})(  p_{tt} , p_t)+(\bm K \nabla p_t,\nabla p_t)+ (\frac{\alpha\beta}{\lambda}-b_0)(  T_{tt} , p_t)\\
& -\frac\beta\lambda(  \xi_{tt} , T_t)+(\frac{\alpha\beta}{\lambda}-b_0)(  p_{tt} , T_t)+(a_0+\frac{\beta^2}{\lambda})(  T_{tt} , T_t)   
  +(\bm\Theta\nabla T_t,\nabla T_t)\\
&=   0.\\
\end{aligned}
\end{equation*}
By the \textbf{Assumption} (A1), the Cauchy-Schwarz inequality and Young's inequality, we have
\begin{equation}\label{norm inequality km thetam C2}
\begin{aligned} 
&  \mu\frac{d}{dt}\Vert\bm\varepsilon( \bm u_t)\Vert_{L^2(\Omega)}^2  
  +\frac1{2 \lambda}\frac{d}{dt}\Vert\alpha  p_t+\beta  T_t- \xi_t\Vert_{L^2(\Omega)}^2
  +\frac{c_0-b_0}{2}\frac{d}{dt} \Vert p_t\Vert_{L^2(\Omega)}^2  
 \\
&+\frac{a_0-b_0}{2}\frac{d}{dt} \Vert T_t\Vert_{L^2(\Omega)}^2 + \frac{b_0}{2} \frac{d}{dt} \Vert  p_t- T_t\Vert_{L^2(\Omega)}^2  
+ k_m\Vert\nabla  p_t\Vert_{L^2(\Omega)}^2
 \\
&+\theta_m \Vert\nabla T_t\Vert_{L^2(\Omega)}^2\le   0.\\
\end{aligned}
\end{equation}
 Integrating equation \eqref{norm inequality km thetam C2} from 0 to $t$, we have
\begin{equation*}    
\begin{aligned} 
&   \mu\Vert\bm\varepsilon( \bm u_t)\Vert_{L^2(\Omega)}^2  
  +\frac1{ 2\lambda}\Vert\alpha  p_t+\beta  T_t- \xi_t\Vert_{L^2(\Omega)}^2
  +\frac{(c_0-b_0)}2 \Vert p_t\Vert_{L^2(\Omega)}^2  
 +\frac{(a_0-b_0)}2 \Vert T_t\Vert_{L^2(\Omega)}^2 \\
&+ \frac{b_0}2  \Vert  p_t- T_t\Vert_{L^2(\Omega)}^2  
+ \int_0^tk_m \Vert\nabla  p_t\Vert_{L^2(\Omega)}^2
 +\int_0^t\theta_m \Vert\nabla  T_t\Vert_{L^2(\Omega)}^2\\
&\le  \mu \Vert\bm\varepsilon( \bm u_t(0))\Vert_{L^2(\Omega)}^2  
  +\frac1{2 \lambda} \Vert\alpha  p_t(0)+\beta  T_t(0)- \xi_t(0)\Vert_{L^2(\Omega)}^2
  +\frac{c_0-b_0}{2}  \Vert p_t(0)\Vert_{L^2(\Omega)}^2  \\
& +\frac{a_0-b_0}{2} \Vert T_t(0)\Vert_{L^2(\Omega)}^2 
+ \frac{b_0}{2}   \Vert  p_t(0)- T_t(0)\Vert_{L^2(\Omega)}^2  ,\\
\end{aligned}
\end{equation*}
which implies that \eqref{ut2+alpha beta C2} holds.
 \end{proof}

With the assistance of Theorem \ref{Th: boundness1}, Lemma \ref{le: energy law}, and Lemma \ref{Th: boundness2}, we can establish the existence and uniqueness of the weak solution to the problem \eqref{TP_Model_four}.
Since the system is linear, the proof can be easily obtained by employing the Galerkin method and the compactness argument, along with the insights presented in \cite{brun2019well}.
\begin{theorem}\label{Th: exist and uniq}
Assume that the {\bf Assumptions} of (A1)-(A3) hold. 
  Let $\bm u^0\in H^1(\Omega)$, $ p^0 \in L^2 (\Omega)$, $T^0\in L^2 (\Omega)$, $\bm f\in L^2 (\Omega)$, $g \in L^2 (\Omega)$, and $H_{s} \in L^2 (\Omega)$. Then a unique weak solution exists to the problem \eqref{TP_Model_four}. Furthermore, if $\bm u^0\in H^2(\Omega)$, $ p^0 \in H^2 (\Omega)$, $T^0\in H^2 (\Omega)$, $\bm f\in L^2 (\Omega)$, $g \in L^2 (\Omega)$, and $H_{s} \in L^2 (\Omega)$, then the energy estimations \eqref{ut2+alpha beta C1} and \eqref{ut2+alpha beta C2}  for the solution hold true.
\end{theorem}

\section{Numerical algorithms  }\label{sec: algorithms}

This section introduces numerical algorithms derived from the four-field formulation previously discussed. Let \( h \) denote the maximum diameter of all elements in the mesh, and for \( h > 0 \), let \(\mathcal{T}_h\) represent a family of triangulations of the domain \(\Omega\) composed of triangular elements. The triangulations are assumed to be shape-regular and quasi-uniform. The finite element spaces defined on \(\mathcal{T}_h\) are as follows:
\[
\begin{aligned}
&\bm{V}_h = \left\{ \bm{v}_h \in [H_{0,\Gamma_d}^1(\Omega)]^2 \cap \bm{C}^0(\bar{\Omega}) : \bm{v}_h|_K \in \mathbb{P}_2(K), \, \forall K \in \mathcal{T}_h \right\}, \\
&Q_h = \left\{ \phi_h \in L^2(\Omega) \cap C^0(\bar{\Omega}) : \phi_h|_K \in \mathbb{P}_1(K), \, \forall K \in \mathcal{T}_h \right\}, \\
&W_h = \left\{ q_h \in H_0^1(\Omega) \cap C^0(\bar{\Omega}) : q_h|_K \in \mathbb{P}_1(K), \, \forall K \in \mathcal{T}_h \right\}.
\end{aligned}
\]
The discrete spaces for \(\bm{u}\) and \(\xi\) are denoted by \(\bm{V}_h \times Q_h\), which are required to satisfy the following fundamental stability condition.
 \begin{equation}\label{infsup for dis}
\begin{aligned} 
 \inf_{\xi\in Q_h}\sup_{\bm v\in \bm{V}_h}\frac{(\nabla\cdot\bm v,\xi)}{\|\bm v\|_{H^1(\Omega)}\|\xi\|_{L^2(\Omega)}}\ge \gamma>0,
\end{aligned}
\end{equation}
 where $\gamma$ is independent of $h$, which means $V_h\times Q_h$ is a stable Stokes pair.  As an example, we utilize the Taylor-Hood elements for the pair \((\bm{u}, \xi)\), specifically \((\bm{P}_2, P_1)\) Lagrange finite elements, and \(P_1\) Lagrange finite elements for both the fluid pressure \(p\) and temperature \(T\).
It is worth noting that other stable Stokes element pairs, such as the MINI element, could also be used. 
 
For time discretization, we consider an equidistant partition consisting of a series of points from $0$ to $T$ with a constant step size of $\Delta t$. Specifically, we define $\bm u^n$ as $\bm u(t_n)$, where $t_n$ represents the $n$-th point in the partition. Similarly, we define $\xi^n=\xi(t_n)$ , $p^n=p(t_n)$ and $T^n=T(t_n)$.
 
\subsection{A coupled algorithm}\label{sub:A coupled algorithm}
Suppose that an initial value $(\bm u_h^0,~\xi_h^0,~p_h^0,~T_h^0)\in \bm V_h\times Q_h\times W_h\times W_h $ is provided.
 Using  the back Euler method to \eqref{TP_Model_var}, for the given finite element spaces $\bm V_h \subset \bm V$, $Q_h\subset Q$, $W_h\subset W $, the coupled algorithm ({\bf Alg. 1}) reads as, find $(\bm u_h^n,~\xi_h^n,~p_h^n,~T_h^n)\in\bm V_h\times Q_h\times W_h\times W_h$
such that for all $(\bm v_h ,~\phi_h ,~q_h ,~S_h )\in\bm V_h\times Q_h\times W_h\times W_h$,
\begin{subequations}
\begin{equation}\label{TP_Model_dis_a}
\begin{aligned} 
2\mu( \bm\varepsilon(\bm u^n_h),\bm\varepsilon(\bm v_h)) - (\nabla\cdot\bm v_h,\xi^n_h)&=(\bm f,\bm v_h)  
\end{aligned}
\end{equation}
\begin{equation}\label{TP_Model_dis_b}
\begin{aligned} 
  -(\nabla\cdot\bm u^n_h,\phi_h)-\frac{1}{\lambda}(\xi^n_h,\phi_h)+\frac{\alpha}{\lambda}(p^n_h,\phi_h)+\frac{\beta}{\lambda}(T^n_h,\phi_h)&=0 ,\\
\end{aligned}
\end{equation}
\begin{equation}\label{TP_Model_dis_c}
\begin{aligned} 
 -  \frac\alpha\lambda(\xi_h^n,q_h)+(c_0+\frac{\alpha^2}{\lambda})(p_h^n,~q_h)
 + (\frac{\alpha\beta}{\lambda}-b_0)(T_h^n,~q_h) 
                +   \Delta t( \bm K\nabla  p^n_h,~\nabla q_h)&= \\
          (c_0+\frac{\alpha^2}{\lambda})(p_h^{n-1},~q_h)
           + (\frac{\alpha\beta}{\lambda}-b_0)(T_h^{n-1},~q_h)
            -  \frac\alpha\lambda(\xi_h^{n-1},~q_h)
           &+\Delta t(g,~q_h),\\
\end{aligned}
\end{equation}
\begin{equation}\label{TP_Model_dis_d}
\begin{aligned}   
  -\frac\beta\lambda(\xi_h^n,S_h)+(a_0+\frac{\beta^2}{\lambda})(T_h^n,S_h)
  +(\frac{\alpha\beta}{\lambda}-b_0)(p_h^n,~S_h)  
            +\Delta t(\bm\Theta\nabla T^n_h,~\nabla S_h) &    =         \\
        (a_0+\frac{\beta^2}{\lambda})(T_h^{n-1},~S_h) 
           +(\frac{\alpha\beta}{\lambda}-b_0)(p_h^{n-1},~S_h)    
               -\frac\beta\lambda(\xi_h^{n-1},~S_h)
           & +\Delta t( H_{s},~S_h).\\
\end{aligned}
\end{equation}
\end{subequations}

To avoid solving the above problem in a fully coupled manner, the system \eqref{TP_Model_dis_a}-\eqref{TP_Model_dis_d} can be reformulated into two decoupled subproblems. If the term \(\frac{\alpha}{\lambda}(p^n_h, \phi_h) + \frac{\beta}{\lambda}(T^n_h, \phi_h)\) in \eqref{TP_Model_dis_b} were moved to the right-hand side, equations \eqref{TP_Model_dis_a} and \eqref{TP_Model_dis_b} become equivalent to the variational form of the generalized Stokes problem for the variables \(\bm{u}_h^n\) and \(\xi_h^n\). Similarly, moving the terms \(\frac{\alpha}{\lambda}(\xi_h^n, q_h)\) in \eqref{TP_Model_dis_c} and \(\frac{\beta}{\lambda}(\xi_h^n, S_h)\) in \eqref{TP_Model_dis_d} to the right-hand sides yields variational forms corresponding to reaction-diffusion equations for \(p_h^n\) and \(T_h^n\) respectively.
These observations facilitate the design of two time-extrapolation-based decoupled algorithms: The first algorithm solves for \((\bm{u}_h^n, \xi_h^n)\) first, followed by solving the reaction-diffusion equations for \((p_h^n, T_h^n)\). Conversely, the second algorithm begins with \((p_h^n, T_h^n)\) and subsequently computes \((\bm{u}_h^n, \xi_h^n)\).

Based on prior experience with Biot’s model \cite{feng2017analysis, ju2020parameter}, these decoupling strategies typically impose stability constraints that necessitate small time step sizes. Specifically, the time step \(\Delta t\) must generally be of order \(O(h^2)\). Consequently, time-extrapolation-based decoupled algorithms may fail to ensure stability or accuracy if the time step size is too large.
    

\subsection{An iteratively decoupled algorithm  }

To overcome the stability constraints, we propose an iterative decoupling algorithm for solving \((\bm{u}_h^n, \xi_h^n, p_h^n, T_h^n)\). We have the following \textbf{Iteratively Decoupled Algorithm (Algorithm 2)}:

For each time step \(t_n\) (\(n \geq 1\)), let the previous step states \((p_h^{n,i-1}, T_h^{n,i-1}, \xi_h^{n,i-1})\) be given, with the initial values set as \(p_h^{n,0} = p_h^{n-1}\), \(T_h^{n,0} = T_h^{n-1}\), and \(\xi_h^{n,0} = \xi_h^{n-1}\). Let \(i\) be the iteration index. For a given \(n\), the \(i\)-th iteration is formulated as:   
 
\textbf{step 1:} Given $\xi_h^{n,i-1}$, find $(p_h^{n,i},T_h^{n,i})$ such that 
\begin{equation}\label{TP_Model_Ts_N1} 
  \begin{aligned}  
       (c_0  + \frac{\alpha^2}{\lambda})(p_{h}^{n,i},q_h) 
            &   + (\frac{\alpha\beta}{\lambda}-b_0)(T_{h}^{n,i},q_h) 
                + \Delta t(\bm K\nabla p^{n,i}_h,\nabla q_h)
                  \\
                =&   \Delta t(g,q_h)  
     +(c_0+\frac{\alpha^2}{\lambda})(p_{h}^{n-1},q_h) 
                + (\frac{\alpha\beta}{\lambda}-b_0)(T_{h}^{n-1},q_h) \\
            &   +  \frac\alpha\lambda(\xi_{h}^{n,i-1},q_h)
                -  \frac\alpha\lambda(\xi_{h}^{n-1},q_h)         ,\\
      (a_0    + \frac{\beta^2}{\lambda})(T_h^{n,i},S_h)
            & +(\frac{\alpha\beta}{\lambda}-b_0)(p_h^{n,i},S_h) 
               +\Delta t(\bm\Theta\nabla T^{n,i}_h,\nabla S_h)\\ 
               =&  \Delta t( H_{s},S_h)                            
    +(a_0    +\frac{\beta^2}{\lambda})(T_h^{n-1},S_h)
               +(\frac{\alpha\beta}{\lambda}-b_0)(p_h^{n-1},S_h)  \\
            & +\frac\beta\lambda(\xi_h^{n,i-1},S_h)
               -  \frac\beta\lambda(\xi_{h}^{n-1},S_h),\\
\end{aligned}  
\end{equation}
  
 \textbf{step 2:} Given $p_h^{n,i},T_h^{n,i}$,find $(\bm u_h^{n,i},\xi_h^{n,i})$ such that 
\begin{equation}\label{TP_Model_uxi_N1} 
\begin{aligned}  
    2\mu(\bm\varepsilon(\bm u^{n,i}_h), \bm\varepsilon(\bm v_h))
  - (\nabla\cdot\bm v_h,\xi^{n,i}_h) 
                &=(\bm f,\bm v_h),\\
   (\nabla\cdot\bm u^{n,i}_h,\phi_h)+\frac{1}{\lambda}(\xi^{n,i}_h,\phi_h)   
                &=  \frac{\alpha}{\lambda}(p^{n,i}_h,\phi_h)+\frac{\beta}{\lambda}(T^{n,i}_h,\phi_h). \\ 
\end{aligned}
\end{equation}

The algorithm proceeds as follows:
1. Compute iteration Step 1 and Step 2 in succession.
2. Repeat the iteration process until the sequence \((\bm{u}_h^{n,i}, \xi_h^{n,i}, p_h^{n,i}, T_h^{n,i})\) converges to within a prescribed tolerance.


For simplicity, the backward Euler scheme has been employed for time discretization in both \textbf{Alg. 1} and \textbf{Alg. 2}. However, it is important to note that alternative higher-order time-stepping methods could also be applied effectively in this context.

\section{Convergence analysis}\label{sec: convergence analysis}
This section establishes the convergence of \textbf{Alg. 2}, demonstrating that the iterative sequences \((\bm u_h^{n,i}, \xi_h^{n,i}, p_h^{n,i}, T_h^{n,i})\) converge to the solution \((\bm u_h^{n}, \xi_h^{n}, p_h^{n}, T_h^{n})\) of the coupled algorithm \textbf{Alg. 1} as \(i \to \infty\). Additionally, we present the error analysis for \textbf{Alg. 1}, confirming its unconditional stability and convergence. Assuming the discrete spaces \(\bm V_h\), \(Q_h\), and \(W_h\), we establish that the time discretization error is of order \(O(\Delta t)\), while the energy-norm errors for \(\bm u\) and \(\xi\) are of order \(O(h^2)\). To begin, we review an essential lemma that aids us in proving the theorem.
\begin{lemma}\label{inverse inequ}
There exists a constant \(C_{inv} > 0\) such that the following inverse inequality is satisfied:  
\[
\|\bm v\|_{H^1(\Omega)} \leq C_{inv} h^{-1} \|\bm v\|_{L^2(\Omega)}, \quad \forall \bm v \in \bm V_h,
\]
where \(h = \max_{K \in \mathcal{T}_h} h_K\), as established in \cite{ciarlet2002finite}.
\end{lemma}

To proceed the analysis, in the remaining part of the article, we
let $e_{\bm u}^i,e_{\xi}^i,e_{p}^i,e_{T}^i$ be differences between the solutions at the iteration $i$ of problem \eqref{TP_Model_Ts_N1}-\eqref{TP_Model_uxi_N1} and
the solutions to problem \eqref{TP_Model_dis_a}-\eqref{TP_Model_dis_d}, i.e.
\begin{equation*}\label{diff_iter_real}
\begin{aligned} 
     (e_{\bm u}^i,e_{\xi}^i,e_{p}^i,e_{T}^i):=(\bm u_h^{n},\xi_h^{n},p_h^{n},T_h^{n})
                                             -(\bm u_h^{n,i},\xi_h^{n,i},p_h^{n,i},T_h^{n,i}).
\end{aligned}
\end{equation*}

Now, we present the following key result, establishing the convergence of {\bf Alg.2}.

\begin{theorem}\label{th: convergence_of_N1}
Let $(\bm u_h^n,\xi_h^n,p_h^n,T_h^n)$ be the solution of the coupled problem \eqref{TP_Model_dis_a}-\eqref{TP_Model_dis_d}. 
Let $(\bm u_h^{n,i},\xi_h^{n,i},p_h^{n,i},T_h^{n,i})$ be the solution of the iteratively decoupled problem \eqref{TP_Model_Ts_N1}-\eqref{TP_Model_uxi_N1}. 
Assuming the \textbf{Assumptions} (A1)-(A3) hold, then, 
the solution sequence $\{(\bm u_h^{n,i},\xi_h^{n,i},p_h^{n,i},T_h^{n,i})\}_{i=1}^{\infty}$ generated by the iteratively decoupled algorithm converges to the solution $(\bm u_h^n,\xi_h^n,p_h^n,T_h^n)$ of the coupled problem. 
\end{theorem}

\begin{proof}

Subtracting equation \eqref{TP_Model_Ts_N1}-\eqref{TP_Model_uxi_N1} from equation \eqref{TP_Model_dis_a}-\eqref{TP_Model_dis_d} respectively, 
we obtain
\begin{subequations}
\begin{equation}\label{eq: var_N1a}
\begin{aligned} 
       (c_0  + \frac{\alpha^2}{\lambda})(e_{p}^i,q_h) 
                + (\frac{\alpha\beta}{\lambda}-b_0)(e_{T}^i,q_h) 
                &+ \Delta t(\bm K\nabla e_{p}^i,\nabla q_h)                        
                =     \frac\alpha\lambda(e_{\xi}^{i-1},q_h)       ,\\
\end{aligned}
\end{equation}
\begin{equation}\label{eq: var_N1b}
\begin{aligned} 
      (a_0    + \frac{\beta^2}{\lambda})(e_{T}^i,S_h)
              +(\frac{\alpha\beta}{\lambda}-b_0)(e_{p}^i,S_h)    
               +\Delta t(\bm\Theta\nabla e_{T}^i,\nabla S_h)        
               =  &   \frac\beta\lambda(e_{\xi}^{i-1},S_h) ,\\
\end{aligned}
\end{equation} 
\begin{equation}\label{eq: var_N1c}
\begin{aligned} 
2\mu(\bm\varepsilon(e_{\bm u}^i),\bm\varepsilon(\bm v_h)) - (\nabla\cdot\bm v_h,e_{\xi}^i)=&0 ,   \\
\end{aligned}
\end{equation}
\begin{equation}\label{eq: var_N1d}
\begin{aligned} 
  (\nabla\cdot e_{\bm u}^{i},\phi_h)+\frac{1}{\lambda}(e_{\xi}^i,\phi_h) 
   =& 
      \frac{\alpha}{\lambda}(e_{p}^i,\phi_h)+\frac{\beta}{\lambda}(e_{T}^i,\phi_h).\\ 
\end{aligned}
\end{equation} 
\end{subequations}
Taking $q_h=e_{p}^i$, $S_h=e_{T}^i$, $\bm v_h=e_{\bm u}^{i}$ and $\phi_h=e_{\xi}^i$ in equation \eqref{eq: var_N1a},  \eqref{eq: var_N1b}, \eqref{eq: var_N1c}  and \eqref{eq: var_N1d}, we have   
\begin{equation}\label{test_q_p_N1a} 
\begin{aligned}  (c_0  +\frac{\alpha^2}{\lambda})\Vert e_{p}^i\|_{L^2(\Omega)}^2
                + (\frac{\alpha\beta}{\lambda}-b_0)( e_{T}^i, e_{p}^i)
                 &+ \Delta t(\bm K\nabla e_{p}^i,\nabla e_{p}^i)    
                =  \frac\alpha\lambda( e_{\xi}^{i-1},e_{p}^i  )     ,\\ 
\end{aligned}
\end{equation}
\begin{equation}\label{test_s_T_N1b}  
\begin{aligned} (a_0    +\frac{\beta^2}{\lambda})\Vert e_{T}^i\|_{L^2(\Omega)}^2
                + (\frac{\alpha\beta}{\lambda}-b_0)( e_{T}^i, e_{p}^i)  
                 + \Delta t(\bm\Theta\nabla e_{T}^i,\nabla e_{T}^i )  
                =&   \frac\beta\lambda( e_{\xi}^{i-1}, e_{T}^i)     .\\ 
\end{aligned}
\end{equation}
\begin{equation}\label{test_v_u_N1}
\begin{aligned} 
2\mu(\bm\varepsilon(e_{\bm u}^i),\bm\varepsilon(e_{\bm u}^{i})) - (\nabla\cdot e_{\bm u}^i,e_{\xi}^i)=&0       .  \\
\end{aligned}
\end{equation} 
\begin{equation}\label{test_phi_xi}
\begin{aligned} 
   (\nabla\cdot e_{\bm u}^{i},e_{\xi}^i)+\frac{1}{\lambda}(e_{\xi}^i,e_{\xi}^i)
   =&\frac{\alpha}{\lambda}(e_{p}^i,e_{\xi}^i)+\frac{\beta}{\lambda}(e_{T}^i,e_{\xi}^i).\\ 
\end{aligned}
\end{equation}  
Summing up equation \eqref{test_q_p_N1a} and equation \eqref{test_s_T_N1b}, we have 
\begin{equation*}\label{eq: sum_N1}  
\begin{aligned}   
     & (c_0 + \frac{\alpha^2}{\lambda})\Vert e_{p}^i\|_{L^2(\Omega)}^2 
     + (a_0 + \frac{\beta^2}{\lambda})\Vert e_{T}^i\|_{L^2(\Omega)}^2 
     + 2(\frac{\alpha\beta}{\lambda}-b_0)(e_{p}^i,e_{T}^i)\\
     &+ \Delta t(\bm K\nabla e_{p}^i,\nabla  e_{p}^i) 
     + \Delta t(\bm\Theta\nabla e_{T}^i,\nabla e_{T}^i )  
          \\
   =&      \frac\alpha\lambda(e_{p}^i,e_{\xi}^{i-1})
       +   \frac\beta\lambda(e_{T}^i,e_{\xi}^{i-1}) 
       =    ( \frac\alpha\lambda e_{p}^i+\frac\beta\lambda e_{T}^i,e_{\xi}^{i-1})
         .\\
\end{aligned}
\end{equation*}
By the Cauchy-Schwarz inequality and \eqref{eq: sum_N1}, we have 
\begin{equation}\label{eq: sc0b0a0 le exii-1}  
\begin{aligned}
& (c_0-b_0)\|e_{p}^i\|_{L^2(\Omega)}^2+(a_0-b_0)\|e_{T}^i\|_{L^2(\Omega)}^2
  +\lambda\|\frac\alpha\lambda e_{p}^i+\frac\beta\lambda e_{T}^i \|_{L^2(\Omega)}^2\\
\le&    c_0 \|e_{p}^i\|_{L^2(\Omega)}^2-2b_0(e_{p}^i,e_{T}^i)+ a_0\|e_{T}^i\|_{L^2(\Omega)}^2         
 +\frac{\alpha^2}{\lambda}\|e_{p}^i\|_{L^2(\Omega)}^2+2\frac{\alpha\beta}{\lambda}(e_{p}^i,e_{T}^i)
+\frac{\beta^2}{\lambda}\|e_{T}^i\|_{L^2(\Omega)}^2                              \\
   \le & (c_0 + \frac{\alpha^2}{\lambda})\Vert e_{p}^i\|_{L^2(\Omega)}^2 
     + (a_0 + \frac{\beta^2}{\lambda})\Vert e_{T}^i\|_{L^2(\Omega)}^2 
     + 2(\frac{\alpha\beta}{\lambda}-b_0)(e_{p}^i,e_{T}^i)\\
     &+ \Delta t(\bm K\nabla e_{p}^i,\nabla  e_{p}^i) 
     + \Delta t(\bm\Theta\nabla e_{T}^i,\nabla e_{T}^i )  
          \\
   =&      ( \frac\alpha\lambda e_{p}^i+\frac\beta\lambda e_{T}^i,e_{\xi}^{i-1})\\
   \le&   \|\frac\alpha\lambda e_{p}^i+\frac\beta\lambda e_{T}^i \|_{L^2(\Omega)}  \|e_{\xi}^{i-1}\|_{L^2(\Omega)}.
   \\
\end{aligned}
\end{equation} 
By simultaneously dividing both sides of the equation \eqref{eq: sc0b0a0 le exii-1} by $\|\frac\alpha\lambda e_{p}^i+\frac\beta\lambda e_{T}^i \|_{L^2(\Omega)}$ which is positive, otherwise we yields \(  \|e_{p}^i \|_{L^2(\Omega)}=\|e_{T}^i \|_{L^2(\Omega)}\), we obtain 
\begin{equation*} 
\begin{aligned}
   \frac{(c_0-b_0)\|e_{p}^i\|_{L^2(\Omega)}^2+(a_0-b_0)\|e_{T}^i\|_{L^2(\Omega)}^2}{ \|\frac\alpha\lambda e_{p}^i+\frac\beta\lambda e_{T}^i \|_{L^2(\Omega)}^2 }\|\frac\alpha\lambda e_{p}^i+\frac\beta\lambda e_{T}^i \|_{L^2(\Omega)}+\lambda\|\frac\alpha\lambda e_{p}^i+\frac\beta\lambda e_{T}^i \|_{L^2(\Omega)}  
& \le     \|e_{\xi}^{i-1}\|_{L^2(\Omega)}. \\
\end{aligned}
\end{equation*}  
Let 
\begin{equation*} 
\mathbb{D}:=\frac{1}{\lambda^2}\frac{(c_0-b_0)\|e_{p}^i\|_{L^2(\Omega)}^2+(a_0-b_0)\|e_{T}^i\|_{L^2(\Omega)}^2}{ \|\frac\alpha\lambda e_{p}^i+\frac\beta\lambda e_{T}^i \|_{L^2(\Omega)}^2 } 
=\frac{(c_0-b_0)\|e_{p}^i\|_{L^2(\Omega)}^2+(a_0-b_0)\|e_{T}^i\|_{L^2(\Omega)}^2}{ \|\alpha e_{p}^i+\beta e_{T}^i \|_{L^2(\Omega)}^2 },
\end{equation*} 
which is positive and independent of $\lambda$, we can next derive  
\begin{equation}\label{epi+eTi<=exii-1}  
\begin{aligned}
   \lambda^2\mathbb{D}\|\frac\alpha\lambda e_{p}^i+\frac\beta\lambda e_{T}^i \|_{L^2(\Omega)}+\lambda\|\frac\alpha\lambda e_{p}^i+\frac\beta\lambda e_{T}^i \|_{L^2(\Omega)}  
& \le     \|e_{\xi}^{i-1}\|_{L^2(\Omega)} \\
\|\frac\alpha\lambda e_{p}^i+\frac\beta\lambda e_{T}^i \|_{L^2(\Omega)}& \le \frac1{\lambda^2\mathbb{D}+\lambda} \|e_{\xi}^{i-1}\|_{L^2(\Omega)}. 
\end{aligned}
\end{equation}  
Summing up equation \eqref{test_v_u_N1} and equation \eqref{test_phi_xi}, we obtain 
\begin{equation}  \label{eq:u+xi le xii-1}    
\begin{aligned}  
    2\mu (\bm\varepsilon(e_{\bm u}^i),\bm\varepsilon(e_{\bm u}^{i}))
     +\frac1\lambda\|e_\xi^i\|_{L^2(\Omega)}^2 
   =&   ( \frac\alpha\lambda e_{p}^i+\frac\beta\lambda e_{T}^i,e_{\xi}^{i})  .\\
\end{aligned}
\end{equation} 
Dropping the first positive term, and using the conclusion of \eqref{epi+eTi<=exii-1}, we have 
\begin{equation*}      
\begin{aligned}  
      \|e_\xi^i\|_{L^2(\Omega)}^2 
   \le& \lambda( \frac\alpha\lambda e_{p}^i+\frac\beta\lambda e_{T}^i,e_{\xi}^{i})
   \le\lambda\|\frac\alpha\lambda e_{p}^i+\frac\beta\lambda e_{T}^i\|_{L^2(\Omega)}\|e_{\xi}^{i}\|_{L^2(\Omega)}
    ,\\ 
\end{aligned}
\end{equation*} 
and then
\begin{equation}  \label{eq:exii le exii-1}    
\begin{aligned}   
  \|e_\xi^i\|_{L^2(\Omega)}  
   \le&\lambda \| \frac\alpha\lambda e_{p}^i+\frac\beta\lambda e_{T}^i\|_{L^2(\Omega)}
   \le\frac\lambda{\lambda^2\mathbb{D}+\lambda}\|e_{\xi}^{i-1}\|_{L^2(\Omega)}.\\
\end{aligned}
\end{equation} 
Therefore, the convergence of $\xi_h^{n, i}$ is proved.
Applying Lemma \ref{Korn's inequa} to equation \eqref{test_v_u_N1}, we see that
\begin{equation*}      
\begin{aligned}   
  2\mu\|\bm\varepsilon(e_{\bm u}^i)\|_{L^2(\Omega)}^2 &=2\mu(\bm\varepsilon(e_{\bm u}^{i}),\bm\varepsilon(e_{\bm u}^{i})) =  (\nabla\cdot e_{\bm u}^i,e_{\xi}^i) \\
  &\le \|\nabla\cdot e_{\bm u}^i\|_{L^2(\Omega)}\| e_{\xi}^i\|_{L^2(\Omega)} 
   \le  C_K\|\bm\varepsilon(e_{\bm u}^i)\|_{L^2(\Omega)}\| e_{\xi}^i\|_{L^2(\Omega)}  , \\
\end{aligned}
\end{equation*} 
which yields the convergence of $\bm u_h^{n,i}$. 
By neglecting the third term on the left-hand side of \eqref{eq: sc0b0a0 le exii-1} and utilizing the result of \eqref{epi+eTi<=exii-1} on the right-hand side of \eqref{eq: sc0b0a0 le exii-1}, we obtain
\begin{equation*}      
\begin{aligned}
 (c_0-b_0)\|e_{p}^i\|_{L^2(\Omega)}^2+(a_0-b_0)\|e_{T}^i\|_{L^2(\Omega)}^2 
&\le     \|\frac\alpha\lambda e_{p}^i+\frac\beta\lambda e_{T}^i \|_{L^2(\Omega)}  \|e_{\xi}^{i-1}\|_{L^2(\Omega)}\\
& \le \frac1{\lambda^2\mathbb{D}+\lambda}\|e_{\xi}^{i-1}\|_{L^2(\Omega)}^2,  \\
\end{aligned}
\end{equation*} 
 which implies the convergence of $p_h^{n,i}$ and $T_h^{n,i}$.
\end{proof}

\begin{remark} \label{re: a0-b0=c0-b0}
The theorem remains valid in the special cases where either \( c_0 - b_0 = 0 \) or \( a_0 - b_0 = 0 \), provided the other term is nonzero, as the proof remains unaffected under these conditions. 

If both \( c_0 - b_0 = 0 \) and \( a_0 - b_0 = 0 \), we can still demonstrate the convergence of the iterative decoupling algorithm. Assuming \( c_0 - b_0 = a_0 - b_0 = 0 \), it follows that \( \mathbb{D} = 0 \). From inequality \eqref{eq:exii le exii-1}, it is evident that \( \|e_{\xi}^i\|_{L^2(\Omega)} \) forms a monotonically non-increasing sequence bounded below by zero, ensuring its convergence. Let 
\[
\lim_{i \to \infty} \|e_{\xi}^i\|_{L^2(\Omega)} = \rho > 0.
\]
From \eqref{eq:u+xi le xii-1} and \eqref{eq:exii le exii-1}, we derive:
\[
2\mu(\bm{\varepsilon}(e_{\bm{u}}^i), \bm{\varepsilon}(e_{\bm{u}}^i)) 
+ \frac{1}{\lambda} \|e_{\xi}^i\|_{L^2(\Omega)}^2 
\leq \|e_{\xi}^i\|_{L^2(\Omega)} \|e_{\xi}^{i-1}\|_{L^2(\Omega)}.
\]
Given \( \rho > 0 \), this implies \( \lim_{i \to \infty} \|\bm{\varepsilon}(e_{\bm{u}}^i)\|_{L^2(\Omega)} = 0 \) as \( i \to \infty \).
Utilizing the discrete inf-sup condition and referencing \eqref{eq: var_N1c}, we obtain:
\[
\gamma \|e_{\xi}^i\|_{L^2(\Omega)} \leq \sup_{\bm{v}_h \in \bm{V}_h} 
\frac{|(\nabla \cdot \bm{v}_h, e_{\xi}^i)|}{\|\bm{v}_h\|_{H^1(\Omega)}}
= \frac{|2\mu(\bm{\varepsilon}(\bm{v}_h), \bm{\varepsilon}(e_{\bm{u}}^i))|}{\|\bm{v}_h\|_{H^1(\Omega)}}
\lesssim \|\bm{\varepsilon}(e_{\bm{u}}^i)\|_{L^2(\Omega)}.
\]
This leads to \( \rho = 0 \), which contradicts the assumption \( \rho > 0 \). Thus, \( \|e_{\xi}^i\|_{L^2(\Omega)} \to 0 \) as \( i \to \infty \).
Consequently, it follows that \( \|e_{\bm{u}}^i\|_{L^2(\Omega)} \to 0 \), \( \|e_p^i\|_{L^2(\Omega)} \to 0 \), and \( \|e_T^i\|_{L^2(\Omega)} \to 0 \) as \( i \to \infty \).
\end{remark}

Here, we denote $d_t \eta_n := ( \eta^n - \eta^{n-1}  )/\Delta t$,~where $\eta$ can be a vector or a scalar.

\begin{lemma}
Let $\{(\bm u_h^n,\xi_h^n,p_h^n,T_h^n)\}$ be solution of the coupled algorithm,~then the following identity holds:
\begin{equation}\label{J_h^l+S_h^l=J_h^0}
\begin{aligned} 
J_h^l+S_h^l=J_h^0, \text{ for  } l\ge 1,
\end{aligned}
\end{equation} 
where for $l\ge 0$,
\begin{equation*}  
\begin{aligned} 
J_h^l:=&\mu\Vert\bm\varepsilon(\bm u_h^l)\Vert_{L^2(\Omega)}^2+\frac1{2\lambda}\Vert\alpha p_h^l+\beta T_h^l-\xi_h^l\Vert_{L^2(\Omega)}^2
  +\frac{c_0-b_0}{2}\Vert p_h^l\Vert_{L^2(\Omega)}^2+\frac{a_0-b_0}{2}\Vert T_h^l\Vert_{L^2(\Omega)}^2\\
  &+\frac{b_0}{2}\Vert p_h^l-T_h^l\Vert_{L^2(\Omega)}^2-(\bm f,\bm u_h^l),
\end{aligned}
\end{equation*} 
and
\begin{equation*}  
\begin{aligned} 
S_h^l:= &\Delta t\sum_{n=1}^l\Big[\Delta t\Big(\mu\Vert \bm\varepsilon(d_t\bm u_h^l)\Vert_{L^2(\Omega)}^2
+\frac1{2\lambda}\Vert d_t(\alpha p_h^l+\beta T_h^l-\xi_h^l)\Vert_{L^2(\Omega)}^2
  +\frac{c_0-b_0}{2}\Vert d_t p_h^l\Vert_{L^2(\Omega)}^2\\
  &+\frac{a_0-b_0}{2}\Vert d_t T_h^l\Vert_{L^2(\Omega)}^2
  +\frac{b_0}{2}\Vert p_h^l-T_h^l\Vert_{L^2(\Omega)}^2 \Big)    \\
        &+\bm K(\nabla p_h^n,\nabla p_h^n)+\bm \Theta(\nabla T_h^n,\nabla T_h^n) 
        - (g,\nabla p_h^n)-(H_{s},\nabla T_h^n)\Big].
\end{aligned}
\end{equation*} 
\end{lemma}
\begin{proof}
Using operator $d_t$ for the equation \eqref{TP_Model_dis_b}, taking $\bm v_h=d_t\bm u_h^n,\phi_h=-\xi_h^n,q_h=p_h^n$ and $S_h=T_h^n$, we have
\begin{equation}\label{TP_Model_dis dtuhn,xihn,phn,Thn}
\begin{aligned} 
2\mu(\bm\varepsilon(\bm u^n_h),\bm\varepsilon(d_t\bm u_h^n)) - (\nabla\cdot d_t\bm u_h^n,\xi^n_h)&=(\bm f,d_t\bm u_h^n)       ,   \\
   (\nabla\cdot d_t\bm u^n_h,\xi_h^n)+\frac{1}{\lambda}( d_t\xi^n_h,\xi_h^n)-\frac{\alpha}{\lambda}( d_tp^n_h,\xi_h^n)-\frac{\beta}{\lambda}( d_tT^n_h,\xi_h^n)&=0 ,\\
 -  \frac\alpha\lambda(\xi_h^n,p_h^n)+(c_0+\frac{\alpha^2}{\lambda})(p_h^n,p_h^n)
 + (\frac{\alpha\beta}{\lambda}-b_0)(T_h^n,p_h^n) 
                +   \Delta t( \bm K\nabla  p^n_h,\nabla p_h^n)&= \\
          (c_0+\frac{\alpha^2}{\lambda})(p_h^{n-1},p_h^n)
           + (\frac{\alpha\beta}{\lambda}-b_0)(T_h^{n-1},p_h^n)
            -  \frac\alpha\lambda(\xi_h^{n-1},p_h^n)
           &+\Delta t(g,p_h^n),\\
  -\frac\beta\lambda(\xi_h^n,T_h^n)+(a_0+\frac{\beta^2}{\lambda})(T_h^n,T_h^n)
  +(\frac{\alpha\beta}{\lambda}-b_0)(p_h^n,T_h^n)  
            +\Delta t(\bm\Theta\nabla T^n_h,\nabla T_h^n) &    =         \\
        (a_0+\frac{\beta^2}{\lambda})(T_h^{n-1},T_h^n) 
           +(\frac{\alpha\beta}{\lambda}-b_0)(p_h^{n-1},T_h^n)    
               -\frac\beta\lambda(\xi_h^{n-1},T_h^n)
           & +\Delta t( H_{s},T_h^n).\\
\end{aligned}
\end{equation}
Summing up the four equations of system \eqref{TP_Model_dis dtuhn,xihn,phn,Thn} and by the definition of $d_t$,
we have
\begin{equation*}    
\begin{aligned} 
&2\mu(\bm\varepsilon(\bm u_h^n),\bm\varepsilon(d_t\bm u_h^n))+\frac{1}{\lambda}(d_t\xi_h^n,\xi_h^n)-\frac{\alpha}{\lambda}(d_tp_h^n,\xi_h^n)-\frac{\beta}{\lambda}(d_tT_h^n,\xi_h^n)\\
& -\frac\alpha\lambda(d_t\xi_h^n ,p_h^n)+(c_0+\frac{\alpha^2}{\lambda})(d_tp_h^n ,p_h^n)+(\frac{\alpha\beta}{\lambda}-b_0)(d_tT_h^n ,p_h^n)+(\bm K \nabla p_h^n,\nabla p_h^n) \\
& -\frac\beta\lambda(d_t\xi_h^n ,T_h^n)+(a_0+\frac{\beta^2}{\lambda})(d_tT_h^n ,T_h^n) +(\frac{\alpha\beta}{\lambda}-b_0)(d_tp_h^n ,T_h^n)  
  +(\bm\Theta\nabla T_h^n,\nabla T_h^n)\\
&= (\bm f,d_t\bm u_h^n) + (g,p_h^n) + ( H_{s},T_h^n).\\
\end{aligned}
\end{equation*}
Then, using the identity
\begin{equation}\label{2x xt=dt x^2+dtx^2}
 2(\eta^n,d_t \eta^n)=d_t\Vert\eta^n\Vert_{L^2(\Omega)}^2
                       +\Delta t\|d_t\eta^n\|_{L^2(\Omega)}^2,
\end{equation}
we have
\begin{equation*}    
\begin{aligned} 
& \mu d_t\Vert\bm\varepsilon(\bm u_h^n)\Vert_{L^2(\Omega)}^2  
  +\frac1{2\lambda}d_t\Vert\alpha p_h^n+\beta T_h^n-\xi_h^n\Vert_{L^2(\Omega)}^2
  +\frac{c_0-b_0}{2}d_t\Vert p_h^n\Vert_{L^2(\Omega)}^2  
 +\frac{a_0-b_0}{2}d_t\Vert T_h^n\Vert_{L^2(\Omega)}^2 \\
&+\frac{b_0}{2}d_t\Vert p_h^n-T_h^n\Vert_{L^2(\Omega)}^2  
+\Delta t\Big(\mu\Vert d_t\bm\varepsilon(\bm u_h^n)\Vert_{L^2(\Omega)}^2
+\frac1{2\lambda}\Vert d_t(\alpha p_h^n+\beta T_h^n-\xi_h^n)\Vert_{L^2(\Omega)}^2\\
 & +\frac{c_0-b_0}{2}\Vert d_t p_h^n\Vert_{L^2(\Omega)}^2 
  +\frac{a_0-b_0}{2}\Vert d_t T_h^n\Vert_{L^2(\Omega)}^2
  +\frac{b_0}{2}\Vert d_t( p_h^n-T_h^n)\Vert_{L^2(\Omega)}^2 \Big) 
  + (\bm K \nabla p_h^n,\nabla p_h^n) \\
 &+(\bm\Theta\nabla T_h^n,\nabla T_h^n)=   (\bm f,d_t\bm u_h^n) + (g,p_h^n) + ( H_{s},T_h^n).
\end{aligned}
\end{equation*}
Applying the summation operator $\Delta t\sum_{n=1}^l$ to both sides of the above equation, we obtain \eqref{J_h^l+S_h^l=J_h^0}. 
\end{proof}

For the discrete finite element spaces \( \bm{V}_h \), \( Q_h \), and \( W_h \) defined in Section \ref{sub:A coupled algorithm}, we introduce the projection operators \( \Pi^{\bm{V}_h}: \bm{V} \to \bm{V}_h \), \( \Pi^{Q_h}: Q \to Q_h \), \( \Pi^{W_{p,h}}: W \to W_h \), and \( \Pi^{W_{T,h}}: W \to W_h \), which satisfy the following equations: for all 
 $(\bm v_h,\phi_h,p_h,T_h)\in \bm V_h\times Q_h\times W_h\times W_h$,
\begin{equation}\label{Pi V}
   2\mu(\bm\varepsilon(\Pi^{\bm V_h}\bm u ),\bm\varepsilon( \bm v_h ))
   - (\nabla\cdot\bm v_h,\Pi^{Q_h} \xi)=2\mu(\bm\varepsilon( \bm u ),\bm\varepsilon( \bm v_h ))
   - (\nabla\cdot\bm v_h, \xi ), 
\end{equation} 
\begin{equation}\label{Pi Q}
 (\nabla\cdot\Pi^{\bm V_h}\bm u,\phi_h)+\frac1\lambda(\Pi^{Q_h} \xi,\phi_h)
         =(\nabla\cdot \bm u,\phi_h)+\frac1\lambda( \xi,\phi_h)
\end{equation} 
\begin{equation}\label{Pi Wp}
   \bm K(\nabla\Pi^{W_{p,h,}}p,\nabla q_h)=\bm K(\nabla p,\nabla q_h),
\end{equation} 
\begin{equation}\label{Pi WT}
   \bm \Theta(\nabla\Pi^{W_{T,h}}T,\nabla S_h)=\bm\Theta (\nabla T,\nabla S_h). 
\end{equation} 
Furthermore, for the operators $\Pi^{\bm V_h}$,~$\Pi^{Q_h} $,~$\Pi^{W_{p,h,}}$ and $\Pi^{W_{T,h}}$ we have following
 properties. 
\begin{equation}\label{Pi V Q properties}
   \|\nabla(\Pi^{\bm V_h}\bm u-\bm u)\|_{L^2(\Omega)}+\|\Pi^{Q_h} \xi-\xi\|_{L^2(\Omega)}
   \lesssim h^2\|\bm u\|_{H^3(\Omega)}+h^2\|\xi\|_{H^2(\Omega)}, 
\end{equation} 
\begin{equation*}\label{Pi Wp properties}
   \| \Pi^{W_{p,h,}}p-p \|_{L^2(\Omega)}+h\|\nabla(\Pi^{W_{p,h,}}p-p)\|_{L^2(\Omega)}
   \lesssim h^2 \|p\|_{H^2(\Omega)},
\end{equation*} 
and
\begin{equation*}\label{Pi WT properties}
     \| \Pi^{W_{T,h}}T-T \|_{L^2(\Omega)}+h\|\nabla(\Pi^{W_{T,h}}T-T)\|_{L^2(\Omega)}
   \lesssim h^2 \|T\|_{H^2(\Omega)}.
\end{equation*} 
 
For convenience, we introduce the following notations:
\begin{equation}\label{e=eI+eh}
\begin{aligned} 
 &e_{\bm u}^n=\bm u^n-\bm u_h^n=( \bm u^n-\Pi^{\bm V_h}\bm u^n)+(\Pi^{\bm V_h}\bm u^n -\bm u_h^n):=e_{\bm u}^{I,n}+e_{\bm u}^{h,n},\\
 &e_{\xi}^n=\xi^n-\xi_h^n=(\xi^n-\Pi^{Q_h}\xi^n)+(\Pi^{Q_h}\xi^n -\xi_h^n):=e_{\xi}^{I,n}+e_{\xi}^{h,n},\\ 
 &e_{p}^n=p^n-p_h^n=(p^n-\Pi^{W_{p,h}} p^n)+(\Pi^{W_{p,h}} p^n -p_h^n):=e_{p}^{I,n}+e_{p}^{h,n},\\ 
 &e_{T}^n=T^n-T_h^n=(T^n-\Pi^{W_{p,h}} T^n)+(\Pi^{W_{p,h}} T^n -T_h^n):=e_{T}^{I,n}+e_{T}^{h,n}.\\
\end{aligned}
\end{equation} 
Next, we formulate a discrete energy law, which serves as a discretized counterpart to the continuous statement presented in Lemma \ref{le: energy law}.
\begin{lemma}
 Let $\{(\bm u_h^n,\xi_h^n,p_h^n,T_h^n)\}$ be defined by the coupled algorithm, then the following
 identity holds:
 \begin{equation}\label{Ehl+Dt+Dt^2=Eh0+c}
\begin{aligned} 
& E_h^l + \Delta t \sum_{n=1}^l \big(  (\bm K\nabla e_p^{h,n},\nabla  e_p^{h,n})
                              +  (\bm\Theta\nabla e_T^{h,n},\nabla  e_T^{h,n})  \big) \\
 &+\Delta t^2\sum_{n=1}^l\Big(\mu\Vert d_t\bm\varepsilon(e_{\bm u}^{h,n})\Vert_{L^2(\Omega)}^2
+\frac1{2\lambda}\Vert d_t(\alpha  e_p^{h,n}+\beta  e_T^{h,n}- e_\xi^{h,n})\Vert_{L^2(\Omega)}^2
  +\frac{c_0-b_0}{2}\Vert d_t  e_p^{h,n}\Vert_{L^2(\Omega)}^2\\
  &+\frac{a_0-b_0}{2}\Vert d_t  e_T^{h,n}\Vert_{L^2(\Omega)}^2
  +\frac{b_0}{2}\Vert  e_p^{h,n}- e_T^{h,n}\Vert_{L^2(\Omega)}^2 \Big)  \\
 =& E_h^0 +\Delta t \sum_{n=1}^l\Big( (\nabla\cdot (d_t\bm u^n- \bm u_t^n ), e_{\xi}^{h,n})
    +\frac{1}{\lambda}(d_t\xi^n- \xi_t^n, e_{\xi}^{h,n})
   -\frac{\alpha}{\lambda}( d_tp^n- p_t^n, e_{\xi}^{h,n})\\
  & -\frac{\beta}{\lambda}(d_tT^n- T_t^n, e_{\xi}^{h,n})
  -\frac\alpha\lambda(d_t\Pi^{Q_h} \xi^n- \xi_t^n , e_p^{h,n})
  +(c_0+\frac{\alpha^2}{\lambda})(d_t\Pi^{W_{p,h,}}p^n- p_t^n , e_p^{h,n}) \\
&+ (\frac{\alpha\beta}{\lambda}-b_0)(d_t\Pi^{W_{T,h}}T^n - T_t^n, e_p^{h,n}) 
  -\frac\beta\lambda(d_t\Pi^{Q_h} \xi^n- \xi_t^n, e_T^{h,n})\\
&+(\frac{\alpha\beta}{\lambda}-b_0)(d_t\Pi^{W_{p,h,}}p^n- p_t^n, e_T^{h,n})
 +(a_0+\frac{\beta^2}{\lambda})(d_t\Pi^{W_{T,h}}T^n - T_t^n, e_T^{h,n})\Big ),   
   \\                                  
\end{aligned}
\end{equation} 
 
where
\begin{equation*}  
\begin{aligned} 
  E_h^n = &\mu\Vert\bm\varepsilon(e_{\bm u}^{h,n})\Vert_{L^2(\Omega)}^2  
  +\frac1{2\lambda} \Vert\alpha  e_p^{h,n}+\beta  e_T^{h,n}- e_\xi^{h,n}\Vert_{L^2(\Omega)}^2
  +\frac{c_0-b_0}{2} \Vert  e_p^{h,n}\Vert_{L^2(\Omega)}^2 \\ 
 &+\frac{a_0-b_0}{2} \Vert  e_T^{h,n}\Vert_{L^2(\Omega)}^2  
+\frac{b_0}{2} \Vert  e_p^{h,n}- e_T^{h,n}\Vert_{L^2(\Omega)}^2.
\end{aligned}
\end{equation*} 
\end{lemma}

\begin{proof}
Utilizing the first equation of the continuous weak formulation \eqref{TP_Model_var}, the finite element formulations \eqref{TP_Model_dis_a} and \eqref{Pi V}, along with the definitions of \( e_{\bm u}^{h,n} \) and \( e_\xi^{h,n} \) as given in equation \eqref{e=eI+eh}, we derive the following relation:  
\begin{equation}\label{dt-d/dt u}
  2\mu(\bm\varepsilon(e_{\bm u}^{h,n}), \bm\varepsilon(\bm v_h)) - (\nabla\cdot\bm v_h, e_\xi^{h,n}) = 0.
\end{equation}  
Furthermore, by combining the second equation of \eqref{TP_Model_var}, the finite element formulation \eqref{TP_Model_dis_b}, and the projection property \eqref{Pi Q}, we obtain
 \begin{equation}\label{dt-d/dt xi}
\begin{aligned} 
 & (\nabla\cdot d_te_{\bm u}^{h,n},\phi_h)+\frac{1}{\lambda}(d_te_\xi^{h,n},\phi_h)
   -\frac{\alpha}{\lambda}(d_te_p^{h,n},\phi_h)-\frac{\beta}{\lambda}(d_te_T^{h,n},\phi_h)\\
  = & (\nabla\cdot (d_t\bm u^n- \bm u_t^n ),\phi_h)
    +\frac{1}{\lambda}(d_t\xi^n- \xi_t^n,\phi_h)
   -\frac{\alpha}{\lambda}( d_t\Pi^{W_{p,h}}p^n- p_t^n,\phi_h)\\
  & -\frac{\beta}{\lambda}(d_t\Pi^{W_{T,h}}T^n- T_t^n,\phi_h).\\
\end{aligned}
\end{equation} 
Using the last two equations of system \eqref{TP_Model_var} 
and of system \eqref{TP_Model_dis_a}-\eqref{TP_Model_dis_d}, equation \eqref{Pi Wp} and \eqref{Pi WT}, we have
 \begin{equation}\label{eq: dt-d/dt p T 1}
\begin{aligned}  
 &-\frac\alpha\lambda(d_te_\xi^{h,n} ,q_h)+(c_0+\frac{\alpha^2}{\lambda})(d_te_p^{h,n} ,q_h)
 +\bm (\bm K\nabla e_p^{h,n},\nabla q_h)+ (\frac{\alpha\beta}{\lambda}-b_0)(d_te_T^{h,n} ,q_h)\\
  =&-\frac\alpha\lambda(d_t\Pi^{Q_h} \xi^n- \xi_t^n ,q_h)
  +(c_0+\frac{\alpha^2}{\lambda})(d_t\Pi^{W_{p,h,}}p^n- p_t^n ,q_h) \\
&+ (\frac{\alpha\beta}{\lambda}-b_0)(d_t\Pi^{W_{T,h}}T^n - T_t^n,q_h) ,\\
\end{aligned}
\end{equation}
and
 \begin{equation}\label{eq: dt-d/dt p T 2}
\begin{aligned}  
& -\frac\beta\lambda(d_te_\xi^{h,n} ,S_h)+(\frac{\alpha\beta}{\lambda}-b_0)(d_te_p^{h,n} ,S_h)
 +(a_0+\frac{\beta^2}{\lambda})(d_te_T^{h,n},S_h)   
  +(\bm\Theta\nabla e_T^{h,n},\nabla S_h) \\
= &-\frac\beta\lambda(d_t\Pi^{Q_h} \xi^n- \xi_t^n,S_h)
+(\frac{\alpha\beta}{\lambda}-b_0)(d_t\Pi^{W_{p,h,}}p^n- p_t^n,S_h)\\
 &+(a_0+\frac{\beta^2}{\lambda})(d_t\Pi^{W_{T,h}}T^n - T_t^n,S_h).   
   \\
\end{aligned}
\end{equation}
 Setting $\bm v_h=d_te_{\bm u}^{h,n}$ in \eqref{dt-d/dt u}, $\phi_h=e_\xi^{h,n}$ in \eqref{dt-d/dt xi} and
  $q_h=e_p^{h,n}$, $S_h=e_T^{h,n}$ in \eqref{eq: dt-d/dt p T 1}-\eqref{eq: dt-d/dt p T 2}, and adding the resulting equations together, we derive
\begin{equation*}    
\begin{aligned}  
& 2\mu(\bm\varepsilon(e_{\bm u}^{h,n} ),\bm\varepsilon(  d_te_{\bm u}^{h,n} ))
   - (\nabla\cdot d_te_{\bm u}^{h,n},e_\xi^{h,n})\\
 &+ (\nabla\cdot d_te_{\bm u}^{h,n}, e_{\xi}^{h,n})+\frac{1}{\lambda}(d_te_\xi^{h,n}, e_{\xi}^{h,n})
   -\frac{\alpha}{\lambda}(d_te_p^{h,n}, e_{\xi}^{h,n})-\frac{\beta}{\lambda}(d_te_T^{h,n}, e_{\xi}^{h,n})\\
 &-\frac\alpha\lambda(d_te_\xi^{h,n} , e_p^{h,n})+(c_0+\frac{\alpha^2}{\lambda})(d_te_p^{h,n} , e_p^{h,n})
 +  (\bm K\nabla e_p^{h,n},\nabla  e_p^{h,n})+ (\frac{\alpha\beta}{\lambda}-b_0)(d_te_T^{h,n} , e_p^{h,n})\\
 & -\frac\beta\lambda(d_te_\xi^{h,n} , e_T^{h,n})+(\frac{\alpha\beta}{\lambda}-b_0)(d_te_p^{h,n} , e_T^{h,n})
 +(a_0+\frac{\beta^2}{\lambda})(d_te_T^{h,n}, e_T^{h,n})   
  +(\bm\Theta\nabla e_T^{h,n},\nabla  e_T^{h,n}) \\
  =& (\nabla\cdot (d_t\bm u^n- \bm u_t^n ),e_{\xi}^{h,n})
    +\frac{1}{\lambda}(d_t\xi^n- \xi_t^n,e_{\xi}^{h,n})
   -\frac{\alpha}{\lambda}( d_t\Pi^{W_{p,h}}p^n- p_t^n,e_{\xi}^{h,n})\\
  & -\frac{\beta}{\lambda}(d_t\Pi^{W_{T,h}}T^n- T_t^n,e_{\xi}^{h,n})\\
  &-\frac\alpha\lambda(d_t\Pi^{Q_h} \xi^n- \xi_t^n ,e_{p}^{h,n})
  +(c_0+\frac{\alpha^2}{\lambda})(d_t\Pi^{W_{p,h,}}p^n- p_t^n ,e_{p}^{h,n}) \\
&+ (\frac{\alpha\beta}{\lambda}-b_0)(d_t\Pi^{W_{T,h}}T^n - T_t^n,e_{p}^{h,n})\\
  &-\frac\beta\lambda(d_t\Pi^{Q_h} \xi^n- \xi_t^n,e_{T}^{h,n})
+(\frac{\alpha\beta}{\lambda}-b_0)(d_t\Pi^{W_{p,h,}}p^n- p_t^n,e_{T}^{h,n})\\
 &+(a_0+\frac{\beta^2}{\lambda})(d_t\Pi^{W_{T,h}}T^n - T_t^n,e_{T}^{h,n}).   
   \\
\end{aligned}
\end{equation*}
 By the identity \eqref{2x xt=dt x^2+dtx^2}, we have
  \begin{equation*}\label{sum dt-d/dt u xi p T}
\begin{aligned}  
&  d_t\Big(\mu\Vert\bm\varepsilon(e_{\bm u}^{h,n})\Vert_{L^2(\Omega)}^2  
  +\frac1{2\lambda} \Vert\alpha  e_p^{h,n}+\beta  e_T^{h,n}- e_\xi^{h,n}\Vert_{L^2(\Omega)}^2
  +\frac{c_0-b_0}{2} \Vert  e_p^{h,n}\Vert_{L^2(\Omega)}^2  
 +\frac{a_0-b_0}{2} \Vert  e_T^{h,n}\Vert_{L^2(\Omega)}^2 \\
&+\frac{b_0}{2} \Vert  e_p^{h,n}- e_T^{h,n}\Vert_{L^2(\Omega)}^2  \Big)
+ (\bm K \nabla  e_p^{h,n},\nabla  e_p^{h,n})
 +(\bm\Theta\nabla  e_T^{h,n},\nabla  e_T^{h,n})\\
 &+\Delta t\Big(\mu\Vert d_t\bm\varepsilon(e_{\bm u}^{h,n})\Vert_{L^2(\Omega)}^2
+\frac1{2\lambda}\Vert d_t(\alpha  e_p^{h,n}+\beta  e_T^{h,n}- e_\xi^{h,n})\Vert_{L^2(\Omega)}^2
  +\frac{c_0-b_0}{2}\Vert d_t  e_p^{h,n}\Vert_{L^2(\Omega)}^2\\
  &+\frac{a_0-b_0}{2}\Vert d_t  e_T^{h,n}\Vert_{L^2(\Omega)}^2
  +\frac{b_0}{2}\Vert  e_p^{h,n}- e_T^{h,n}\Vert_{L^2(\Omega)}^2 \Big)  \\
  =& (\nabla\cdot (d_t\bm u^n- \bm u_t^n ),e_{\xi}^{h,n})
    +\frac{1}{\lambda}(d_t\xi^n- \xi_t^n,e_{\xi}^{h,n})
   -\frac{\alpha}{\lambda}( d_t\Pi^{W_{p,h}}p^n- p_t^n,e_{\xi}^{h,n})\\
  & -\frac{\beta}{\lambda}(d_t\Pi^{W_{T,h}}T^n- T_t^n,e_{\xi}^{h,n})\\
  &-\frac\alpha\lambda(d_t\Pi^{Q_h} \xi^n- \xi_t^n ,e_{p}^{h,n})
  +(c_0+\frac{\alpha^2}{\lambda})(d_t\Pi^{W_{p,h,}}p^n- p_t^n ,e_{p}^{h,n}) \\
&+ (\frac{\alpha\beta}{\lambda}-b_0)(d_t\Pi^{W_{T,h}}T^n - T_t^n,e_{p}^{h,n})\\
  &-\frac\beta\lambda(d_t\Pi^{Q_h} \xi^n- \xi_t^n,e_{T}^{h,n})
+(\frac{\alpha\beta}{\lambda}-b_0)(d_t\Pi^{W_{p,h,}}p^n- p_t^n,e_{T}^{h,n})\\
 &+(a_0+\frac{\beta^2}{\lambda})(d_t\Pi^{W_{T,h}}T^n - T_t^n,e_{T}^{h,n}).   
   \\
\end{aligned}
\end{equation*}
Then, applying the summation operator $\Delta t\sum_{n=1}^l$ to both sides of the above equation, we get \eqref{Ehl+Dt+Dt^2=Eh0+c}.
\end{proof}

Divide the interval \([0, \tau]\) into \(l\) equal parts and let \(t_n\) represent the time value at the $n$-th time step, \(t_0=0\) and \(t_l=\tau\). 
\begin{theorem}\label{euhn estimate}
 Let $\{(\bm u_h^n, \xi_h^n,p_h^n, T_h^n)\}$ be the solution of  the coupled algorithm, then the following
 estimates hold:
\begin{equation}\label{max Enl + ehn le Dt^2 H^4}
\begin{aligned} 
& \max_{0\le n\le l}\Big[
   \mu\Vert\bm\varepsilon(e_{\bm u}^{h,n})\Vert_{L^2(\Omega)}^2  
  +\frac1{2\lambda} \Vert\alpha  e_p^{h,n}+\beta  e_T^{h,n}- e_\xi^{h,n}\Vert_{L^2(\Omega)}^2
  +\frac{c_0-b_0}{2} \Vert  e_p^{h,n}\Vert_{L^2(\Omega)}^2  \\
 & +\frac{a_0-b_0}{2} \Vert  e_T^{h,n}\Vert_{L^2(\Omega)}^2 
  +\frac{b_0}{2} \Vert  e_p^{h,n}- e_T^{h,n}\Vert_{L^2(\Omega)}^2
                   \Big] \\
                   &+  \Delta t \sum_{n=1}^l \big(  (\bm K\nabla e_p^{h,n},\nabla  e_p^{h,n})
                              +  (\bm\Theta\nabla e_T^{h,n},\nabla  e_T^{h,n})  \big)\\
\le &  \tilde C_1\Delta t^2 + \tilde C_2 h^4,             
\end{aligned}
\end{equation} 
 where
 \begin{equation*}  
\begin{aligned} 
&  \tilde C_1   =\tilde C_1(  \int_{t_0}^{\tau}\| \bm u_{tt}  \|_{H^1(\Omega)}^2,
              \int_{t_0}^{\tau}\| \xi_{tt}  \|_{L^2(\Omega)}^2,
              \int_{t_0}^{\tau}\| p_{tt}  \|_{L^2(\Omega)}^2,
            \int_{t_0}^{\tau}\| T_{tt}  \|_{L^2(\Omega)}^2)\\ 
&  \tilde C_2   =\tilde C_2 ( \int_{t_0}^{\tau}\| \xi_t  \|_{H^2(\Omega)}^2 ,
           \int_{t_0}^{\tau}\| p_t  \|_{H^2(\Omega)}^2, 
             \int_{t_0}^{\tau}\| T_t  \|_{H^2(\Omega)}^2),             
\end{aligned}
\end{equation*} 
\end{theorem}

\begin{proof}
For the initial conditions, we set  
\(\bm u_h^0 = \Pi^{\bm V_h} \bm u^0\), \(\xi_h^0 = \Pi^{Q_h} \xi^0\), \(p_h^0 = \Pi^{W_{p,h}} p^0\), and \(T_h^0 = \Pi^{W_{T,h}} T^0\).  
Using these definitions, the following inequality is obtained from \eqref{Ehl+Dt+Dt^2=Eh0+c}.
\begin{equation}\label{Ehl+Dt+Dt^2=Eh0+c uh0}
\begin{aligned} 
 & E_h^l + \Delta t \sum_{n=1}^l \big(  (\bm K\nabla e_p^{h,n},\nabla  e_p^{h,n})
                              +  (\bm\Theta\nabla e_T^{h,n},\nabla  e_T^{h,n})  \big) \\
 \le&  \Delta t \sum_{n=1}^l\Big( (\nabla\cdot (d_t\bm u^n- \bm u_t^n ), e_{\xi}^{h,n})
    +\frac{1}{\lambda}(d_t\xi^n- \xi_t^n, e_{\xi}^{h,n})
   -\frac{\alpha}{\lambda}( d_tp^n- p_t^n, e_{\xi}^{h,n})\\
  & -\frac{\beta}{\lambda}(d_tT^n- T_t^n, e_{\xi}^{h,n})\\
  &-\frac\alpha\lambda(d_t\Pi^{Q_h} \xi^n- \xi_t^n , e_p^{h,n})
  +(c_0+\frac{\alpha^2}{\lambda})(d_t\Pi^{W_{p,h,}}p^n- p_t^n , e_p^{h,n}) \\
&+ (\frac{\alpha\beta}{\lambda}-b_0)(d_t\Pi^{W_{T,h}}T^n - T_t^n, e_p^{h,n}) ,\\
  &-\frac\beta\lambda(d_t\Pi^{Q_h} \xi^n- \xi_t^n, e_T^{h,n})
+(\frac{\alpha\beta}{\lambda}-b_0)(d_t\Pi^{W_{p,h,}}p^n- p_t^n, e_T^{h,n})\\
 &+(a_0+\frac{\beta^2}{\lambda})(d_t\Pi^{W_{T,h}}T^n - T_t^n, e_T^{h,n}) \Big)\\  
 :=&D_1+D_2 +D_3 +D_4 +D_5 +D_6 +D_7 +D_8 +D_9 +D_{10}.  \\                                  
\end{aligned}
\end{equation} 
We now estimate each term on the right-hand side of \eqref{Ehl+Dt+Dt^2=Eh0+c uh0}.
Using the Cauchy-Schwarz inequality, Taylor series expansion, and Proposition 3.1. in \cite{gu2023priori}, we can bound
$D_1$ by
   \begin{equation*}    
\begin{aligned} 
 D_1= &\Delta t \sum_{n=1}^l   (\nabla\cdot (d_t\bm u^n- \bm u_t^n ), e_{\xi}^{h,n}) \\
 =& \sum_{n=1}^l  (\nabla\cdot (\bm u^n-\bm u^{n-1}-\Delta t \bm u_t^n ), e_{\xi}^{h,n})\\                             
 \lesssim& \frac{1}{\Delta t}   \sum_{n=1}^l 
        \| \bm u^n-\bm u^{n-1}-\Delta t \bm u_t^n  \|_{H^1(\Omega)}^2
             +    \Delta t \sum_{n=1}^l   \| e_{\xi}^{h,n} \|_{L^2(\Omega)}^2\\
 =&   \frac{1}{\Delta t} \sum_{n=1}^l
       \|\int_{t_{n-1}}^{t_n}\bm u_{tt}(s)(t_{n-1}-s)ds  \|_{H^1(\Omega)}^2
       +     \Delta t \sum_{n=1}^l        \| e_{\xi}^{h,n} \|_{L^2(\Omega)}^2  \\
 \lesssim&   
      \Delta t^2 \int_{t_0}^{\tau}\| \bm u_{tt}  \|_{H^1(\Omega)}^2
    + \Delta t \sum_{n=1}^l        \| e_{\xi}^{h,n} \|_{L^2(\Omega)}^2  \\        
\end{aligned}
\end{equation*} 
Similarly,   $D_2$ ,$D_3$ and  $D_4$ can be bounded as follows:  
\begin{equation*}  
   D_2=\Delta t \sum_{n=1}^l \frac{1}{\lambda}(d_t\xi^n- \xi_t^n, e_{\xi}^{h,n}) 
\lesssim   
    \frac1\lambda( \Delta t^2 \int_{t_0}^{\tau}\| \xi_{tt}  \|_{L^2(\Omega)}^2
    + \Delta t \sum_{n=1}^l        \| e_{\xi}^{h,n} \|_{L^2(\Omega)}^2  ),\\                                  
\end{equation*} 
\begin{equation*} 
   D_3=\Delta t \sum_{n=1}^l \frac{\alpha}{\lambda}(d_tp^n- p_t^n, e_{\xi}^{h,n}) 
\lesssim   
     \frac\alpha\lambda( \Delta t^2 \int_{t_0}^{\tau}\| p_{tt}  \|_{L^2(\Omega)}^2
    + \Delta t \sum_{n=1}^l        \| e_{\xi}^{h,n} \|_{L^2(\Omega)}^2 ) ,\\  
\end{equation*} 
and
\begin{equation*}  
   D_4=\Delta t \sum_{n=1}^l \frac{\beta}{\lambda}(d_tT^n- T_t^n, e_{\xi}^{h,n}) 
\lesssim   
    \frac\beta\lambda(  \Delta t^2 \int_{t_0}^{\tau}\| T_{tt}  \|_{L^2(\Omega)}^2
    + \Delta t \sum_{n=1}^l        \| e_{\xi}^{h,n} \|_{L^2(\Omega)}^2 ) .\\                                  
\end{equation*} 
By use of the estimate \eqref{Pi V Q properties}, the Cauchy-Schwarz inequality,
Taylor series expansion, and Proposition 3.1. in \cite{gu2023priori}, we see that $D_5$ satisfies
\begin{equation*}    
\begin{aligned}  
 D_5=& - \Delta t \sum_{n=1}^l\frac\alpha\lambda(d_t\Pi^{Q_h} \xi^n- \xi_t^n , e_p^{h,n}) \\
 =& -\sum_{n=1}^l\frac\alpha\lambda( \Pi^{Q_h} \xi^n-\Pi^{Q_h} \xi^{n-1}-\Delta t \xi_t^n , e_p^{h,n})\\
 =&-\sum_{n=1}^l\frac\alpha\lambda( \Pi^{Q_h} \xi^n-\Pi^{Q_h} \xi^{n-1}-(\xi^n-\xi^{n-1})
         +(\xi^n-\xi^{n-1})-\Delta t \xi_t^n , e_p^{h,n})\\
  \lesssim&
   \frac\alpha\lambda(h^4
     \frac{1}{\Delta t} \sum_{n=1}^l     \| \xi^n-\xi^{n-1}  \|_{H^2(\Omega)}^2
   + \Delta t^2  \int_{t_0}^{\tau}\| \xi_{tt}  \|_{L^2(\Omega)}^2
    + \Delta t \sum_{n=1}^l        \| e_{p}^{h,n} \|_{L^2(\Omega)}^2
   )    \\    
  \lesssim&
   \frac\alpha\lambda(h^4
    \int_{t_0}^{\tau}\| \xi_t  \|_{H^2(\Omega)}^2
   + \Delta t^2  \int_{t_0}^{\tau}\| \xi_{tt}  \|_{L^2(\Omega)}^2
    + \Delta t \sum_{n=1}^l        \| e_{p}^{h,n} \|_{L^2(\Omega)}^2
   ),                           
\end{aligned}
\end{equation*} 
and similarly,
  \begin{equation*}    
\begin{aligned}  
 D_6= & 
   (c_0+\frac{\alpha^2}{\lambda})(d_t\Pi^{W_{p,h,}}p^n- p_t^n , e_p^{h,n})    \\
   \lesssim&
    (c_0+\frac{\alpha^2}{\lambda})(h^4
    \int_{t_0}^{\tau}\| p_t  \|_{H^2(\Omega)}^2
   + \Delta t^2  \int_{t_0}^{\tau}\| p_{tt}  \|_{L^2(\Omega)}^2
    + \Delta t \sum_{n=1}^l        \| e_{p}^{h,n} \|_{L^2(\Omega)}^2
   ), \\                                  
\end{aligned}
\end{equation*} 

  \begin{equation*}    
\begin{aligned}  
 D_7=&  (\frac{\alpha\beta}{\lambda}-b_0)(d_t\Pi^{W_{T,h}}T^n - T_t^n, e_p^{h,n})  \\
 \lesssim&
   (\frac{\alpha\beta}{\lambda}-b_0)(h^4
    \int_{t_0}^{\tau}\| T_t  \|_{H^2(\Omega)}^2
   + \Delta t^2  \int_{t_0}^{\tau}\| T_{tt}  \|_{L^2(\Omega)}^2
    + \Delta t \sum_{n=1}^l        \| e_{p}^{h,n} \|_{L^2(\Omega)}^2
   ),
\end{aligned}
\end{equation*} 

  \begin{equation*}    
\begin{aligned}  
 D_8=& -\frac\beta\lambda(d_t\Pi^{Q_h} \xi^n- \xi_t^n, e_T^{h,n}) \\
     \lesssim&
  \frac\beta\lambda(h^4
    \int_{t_0}^{\tau}\| \xi_t  \|_{H^2(\Omega)}^2
   + \Delta t^2  \int_{t_0}^{\tau}\| \xi_{tt}  \|_{L^2(\Omega)}^2
    + \Delta t \sum_{n=1}^l        \| e_{T}^{h,n} \|_{L^2(\Omega)}^2
   ),                    
\end{aligned}
\end{equation*} 
  \begin{equation*}    
\begin{aligned}  
 D_9=&  (\frac{\alpha\beta}{\lambda}-b_0)(d_t\Pi^{W_{p,h,}}p^n- p_t^n, e_T^{h,n})\\ 
  \lesssim&
   |\frac{\alpha\beta}{\lambda}-b_0|(h^4
    \int_{t_0}^{\tau}\| p_t  \|_{H^2(\Omega)}^2
   + \Delta t^2  \int_{t_0}^{\tau}\| p_{tt}  \|_{L^2(\Omega)}^2
    + \Delta t \sum_{n=1}^l        \| e_{T}^{h,n} \|_{L^2(\Omega)}^2
   )                                   
\end{aligned}
\end{equation*} 

  \begin{equation*}    
\begin{aligned}  
 D_{10}= & (a_0+\frac{\beta^2}{\lambda})(d_t\Pi^{W_{T,h}}T^n - T_t^n, e_T^{h,n}) \\
   \lesssim&(a_0+\frac{\beta^2}{\lambda})(h^4
    \int_{t_0}^{\tau}\|T_t  \|_{H^2(\Omega)}^2
   + \Delta t^2  \int_{t_0}^{\tau}\| T_{tt}  \|_{L^2(\Omega)}^2
    + \Delta t \sum_{n=1}^l        \| e_{T}^{h,n} \|_{L^2(\Omega)}^2
   ).                               
\end{aligned}
\end{equation*} 
The above bounds by the discrete Gronwall’s inequality lead to
  \begin{equation*} 
\begin{aligned} 
 & E_h^l + \Delta t \sum_{n=1}^l \big(  (\bm K\nabla e_p^{h,n},\nabla  e_p^{h,n})
                              +  (\bm\Theta\nabla e_T^{h,n},\nabla  e_T^{h,n})  \big) \\
 \lesssim &\Delta t^2 \int_{t_0}^{\tau}\| \bm u_{tt}  \|_{H^1(\Omega)}^2
           +\Delta t^2 \int_{t_0}^{\tau}\| \xi_{tt}  \|_{L^2(\Omega)}^2
            +\Delta t^2 \int_{t_0}^{\tau}\| p_{tt}  \|_{L^2(\Omega)}^2
          +\Delta t^2 \int_{t_0}^{\tau}\| T_{tt}  \|_{L^2(\Omega)}^2\\
         & +h^4 \int_{t_0}^{\tau}\| \xi_t  \|_{H^2(\Omega)}^2 
         +h^4 \int_{t_0}^{\tau}\| p_t  \|_{H^2(\Omega)}^2 
          +h^4  \int_{t_0}^{\tau}\| T_t  \|_{H^2(\Omega)}^2, \\                                  
\end{aligned}
\end{equation*} 
which implies that \eqref{max Enl + ehn le Dt^2 H^4} holds.
\end{proof}

\begin{theorem}
 Let $\{(\bm u_h^n,\xi_h^n,p_h^n,T_h^n)\}$ be defined by the coupled algorithm,then the following
 estimates hold:
\begin{equation*}\label{max Enl + en le Dt^2 H^4}
\begin{aligned} 
& \max_{0\le n\le l}\Big[
   \mu\Vert\bm\varepsilon(e_{\bm u}^n)\Vert_{L^2(\Omega)}^2  
  +\frac1{2\lambda} \Vert e_\xi^n\Vert_{L^2(\Omega)}^2
  +\frac{c_0-b_0}{2} \Vert  e_p^n\Vert_{L^2(\Omega)}^2  
  +\frac{a_0-b_0}{2} \Vert  e_T^n\Vert_{L^2(\Omega)}^2 
                   \Big]\\
\le &  \tilde C_1\Delta t^2 + \tilde C_2 h^4,             
\end{aligned}
\end{equation*} 
 where $\tilde C_1$ and $\tilde C_2$ are defined the same as in Theorem \ref{euhn estimate}.

\end{theorem}
\begin{proof} 
By use of the \eqref{e=eI+eh} and Theorem \ref{euhn estimate}, we obtain the estimate \eqref{max Enl + en le Dt^2 H^4}.

\end{proof}
 
\section{Numerical experiments}\label{sec: experiments}

This section presents numerical experiments to assess the computational accuracy and efficiency of the proposed algorithms in Section \ref{sec: algorithms}. The objective is to evaluate the performance of the two algorithms under different physical parameter settings. A uniform triangular partition, \(\mathcal{T}_h\), is employed to generate the finite element mesh over the domain \(\Omega = [0,1] \times [0,1]\), with the terminal time set to 
 \(\tau = 0.01\) . The time interval is defined as \(J = [0, \tau]\). The experiments are conducted using the following benchmark problem. 
All computations are carried out using the FreeFEM++ software \cite{Hecht2012NewDI}.

{\bfseries Example 1:}
We appropriately select the body force \(\bm f\), mass source \(g\), and heat source \(h\) to ensure that the exact solution is
\begin{equation} \label{Ex: example1}
\begin{aligned}   
      & \bm u(x,y,t)= e^{-t}\left(\begin{array}{c}  
     \sin(2\pi y)(\cos(2\pi x) -1)+\frac{1}{\mu+\lambda}\sin(\pi x) \sin(\pi y)  \\  
     \sin(2\pi x)(1-\cos(2\pi y) )+\frac{1}{\mu+\lambda}\sin(\pi x) \sin(\pi y)   \\
                                 \end{array}\right),   \\
      & p(x,y,t) =   e^{-t}\sin(\pi x) \sin(\pi y) ,      \\
      & T(x,y,t) =   e^{-t}\sin(\pi x) \sin(\pi y) .
\end{aligned}
\end{equation} 
We impose Dirichlet boundary conditions $\bm{u} = \bm{0}$ on the boundary segment $\Gamma_d = \{(x,y)\,|\, 0 \le y \le 1,\, x = 0 \text{ or } 1\}$, and Neumann boundary conditions on $\Gamma_n = \{(x,y)\,|\, 0 \le x \le 1,\, y = 0 \text{ or } 1\}$. For the pressure $p$ and temperature $T$, homogeneous Dirichlet boundary conditions are prescribed on the entire boundary $\partial\Omega$. The initial conditions for $\bm{u}^0$, $p^0$, and $T^0$ at time $t = 0$ are given by the expression in \eqref{Ex: example1}.


In the experiments, we use uniform triangular meshes with sizes \(h = 1/16\), \(h = 1/32\), \(h = 1/64\), and \(h = 1/128\). These meshes are refined by bisecting each triangle through its midpoints. At the terminal time \(\tau\), we report the computed \(L^2\)-norm errors for \(\xi\) and \(H^1\)-norm errors for \(\bm u\), \(p\), and \(T\), along with their corresponding convergence 
orders. The term "iters" denotes the number of iterations used in the iterative decoupling algorithm. Larger time step sizes are intentionally selected in the tests to emphasize the decoupled algorithm's efficiency and robustness while demonstrating its stability in the absence of restrictive constraints.

For the coupled algorithm (\textbf{Alg. 1}), the time step size is set to \(\Delta t = 10^{-3}\), requiring 10 steps to reach the terminal time \(\tau\). This means that the coupled system \eqref{TP_Model_dis_a}–\eqref{TP_Model_dis_d} is solved 10 times. For the decoupled algorithm (\textbf{Alg. 2}), two scenarios are considered: \(\Delta t = 5 \times 10^{-3}\) and \(\Delta t = 10^{-2}\), corresponding to 200 and 100 time steps, respectively, to reach \(\tau\). In these cases, the number of iterations is set to 5 and 10, respectively. Notably, the total number of system solves for \eqref{TP_Model_Ts_N1}–\eqref{TP_Model_uxi_N1} matches that of \textbf{Alg. 1}, ensuring an equivalent computational cost. This balance is achieved by appropriately selecting \(\Delta t\) values and iteration counts. Despite this equivalence in system solves, the computational time for the decoupled algorithm is significantly reduced because solving smaller, independent subproblems is more efficient than solving the coupled system.

From the results, we have the 
following three key observations:
1. Both the coupled and decoupled algorithms attain optimal convergence rates in the energy norm, as demonstrated in Tables~\ref{tab: K Theta 0.1 nu 0.3} to \ref{tab:a0=b0=c0=0}. Notably, the decoupled algorithm achieves comparable accuracy to the coupled counterpart while reducing computational cost. For instance, CPU times (in seconds) for each algorithm are presented in Table~\ref{tab: K Theta 0.1 nu 0.3} for comparison.
2. The final entries in Tables~\ref{tab: K Theta 0.1 nu 0.3} through \ref{tab:a0=b0=c0=0} indicate that increasing the number of iterations in the decoupled algorithm from 5 to 10 results in slightly improved accuracy and enhanced convergence rates. This highlights the potential for accuracy gains with more iterations in each time step.
3. Both algorithms exhibit robust performance even under extreme parameter regimes, as evidenced in the subsequent numerical experiments. 

These results validate the computational efficiency and stability of the proposed algorithms across varying scenarios.

\subsection{Tests for the parameter $\nu$}

\begin{table}[htbp] \tiny  
\centering
\caption{The errors and convergence orders with decreasing mesh sizes. $ \nu=0.3$. }
\begin{tabular}{cccccccccc}
\toprule
  $h$&$\Vert\bm u\Vert_{H^1}$ &orders &$\Vert\xi\Vert_{L^2}$ &orders &$\Vert p\Vert_{H^1}$ &orders&$\Vert T\Vert_{H^1}$ &orders&CPU time(s) \\ 
\midrule
Alg.1,\\
$\Delta t=10^{-3}$ &\\
\midrule
$1/16$ & 1.00607e-01 & 0.00 & 6.00958e-03 & 0.00 & 2.28033e-01 & 0.00 & 2.28033e-01 & 0.00 &3.32\\
$1/32$ & 2.53649e-02 & 1.99 & 1.48475e-03 & 2.02 & 1.09515e-01 & 1.06 & 1.09515e-01 & 1.06 &8.87 \\
$1/64$ & 6.35806e-03 & 2.00 & 3.69965e-04 & 2.00 & 5.41728e-02 & 1.02 & 5.41728e-02 & 1.02 &31.26\\
$1/128$ & 1.59098e-03 & 2.00 & 9.22429e-05 & 2.00 & 2.70126e-02 & 1.00 & 2.70126e-02 & 1.00 &153.86\\ 
\midrule
Alg.2,\\
$\Delta t=5\times 10^{-3}$,\\
iters=5 &\\
\midrule
$1/16$ & 1.00607e-01 & 0.00 & 6.01533e-03 & 0.00 & 2.28993e-01 & 0.00 & 2.28993e-01 & 0.00 &2.03\\
$1/32$ & 2.53650e-02 & 1.99 & 1.48612e-03 & 2.02 & 1.09638e-01 & 1.06 & 1.09638e-01 & 1.06& 6.67\\
$1/64$ & 6.35808e-03 & 2.00 & 3.70229e-04 & 2.01 & 5.41866e-02 & 1.02 & 5.41866e-02 & 1.02&19.26 \\
$1/128$ & 1.59097e-03 & 2.00 & 9.22425e-05 & 2.00 & 2.70136e-02 & 1.00 & 2.70136e-02 & 1.00 &79.37\\
\midrule
Alg.2,\\
$\Delta t=10^{-2}$,\\
iters=10 &\\
\midrule
$1/16$ & 1.00608e-01 & 0.00 & 6.02199e-03 & 0.00 & 2.30330e-01 & 0.00 & 2.30330e-01 & 0.00& 2.27\\
$1/32$ & 2.53650e-02 & 1.99 & 1.48635e-03 & 2.02 & 1.09772e-01 & 1.07 & 1.09772e-01 & 1.07&7.49 \\
$1/64$ & 6.35790e-03 & 2.00 & 3.68857e-04 & 2.01 & 5.41797e-02 & 1.02 & 5.41797e-02 & 1.02& 21.37\\
$1/128$ & 1.59078e-03 & 2.00 & 9.07180e-05 & 2.02 & 2.70029e-02 & 1.00 & 2.70029e-02 & 1.00&85.17 \\
\bottomrule
\label{tab: K Theta 0.1 nu 0.3}
\end{tabular}            
\end{table} 

 \begin{table}\tiny  
 \centering 
\caption{The errors and convergence orders with decreasing mesh sizes. $ \nu=0.49999$. }
\begin{tabular}{ccccccccc}
\toprule
  $h$&$\Vert\bm u\Vert_{H^1}$ &orders &$\Vert\xi\Vert_{L^2}$ &orders &$\Vert p\Vert_{H^1}$ &orders&$\Vert T\Vert_{H^1}$ &orders \\ 
\midrule
Alg.1,\\
$\Delta t=10^{-3}$\\
\midrule
$1/16$ & 9.99038e-02 & 0.00 & 9.71217e-03 & 0.00 & 2.15222e-01 & 0.00 & 2.15222e-01 & 0.00 \\
$1/32$ & 2.51776e-02 & 1.99 & 2.38401e-03 & 2.03 & 1.07872e-01 & 1.00 & 1.07872e-01 & 1.00 \\
$1/64$ & 6.31033e-03 & 2.00 & 5.93579e-04 & 2.01 & 5.39689e-02 & 1.00 & 5.39689e-02 & 1.00 \\
$1/128$ & 1.57899e-03 & 2.00 & 1.48251e-04 & 2.00 & 2.69886e-02 & 1.00 & 2.69886e-02 & 1.00 \\
\midrule
Alg.2,\\
$\Delta t=5\times 10^{-3}$, \\
iters=5 &\\
\midrule
$1/16$ & 9.99038e-02 & 0.00 & 9.71217e-03 & 0.00 & 2.15234e-01 & 0.00 & 2.15234e-01 & 0.00 \\
$1/32$ & 2.51776e-02 & 1.99 & 2.38401e-03 & 2.03 & 1.07874e-01 & 1.00 & 1.07874e-01 & 1.00 \\
$1/64$ & 6.31033e-03 & 2.00 & 5.93579e-04 & 2.01 & 5.39692e-02 & 1.00 & 5.39692e-02 & 1.00 \\
$1/128$ & 1.57899e-03 & 2.00 & 1.48251e-04 & 2.00 & 2.69886e-02 & 1.00 & 2.69886e-02 & 1.00 \\
\midrule
Alg.2,\\
$\Delta t=10^{-2}$,\\
iters=10 &\\
\midrule
$1/16$ & 9.99038e-02 & 0.00 & 9.71217e-03 & 0.00 & 2.15248e-01 & 0.00 & 2.15248e-01 & 0.00 \\
$1/32$ & 2.51776e-02 & 1.99 & 2.38401e-03 & 2.03 & 1.07876e-01 & 1.00 & 1.07876e-01 & 1.00 \\
$1/64$ & 6.31033e-03 & 2.00 & 5.93579e-04 & 2.01 & 5.39694e-02 & 1.00 & 5.39694e-02 & 1.00 \\
$1/128$ & 1.57899e-03 & 2.00 & 1.48251e-04 & 2.00 & 2.69888e-02 & 1.00 & 2.69888e-02 & 1.00 \\
\bottomrule
\label{tab: nu 0.49999}
\end{tabular}                  
\end{table} 
 
In this subsection, we evaluate the performance of the algorithms under various settings of the Poisson ratio \(\nu\). For the other parameters, we fix \(E = 1\), \(c_0 = a_0 = 0.2\), \(b_0 = 0.1\), \(\alpha = \beta = 0.1\), and \(\bm K = \bm \Theta = 0.1 \bm I\), where \(\bm I\) denotes the \(2 \times 2\) identity matrix. The data shown in Table \ref{tab: K Theta 0.1 nu 0.3}, along with the subsequent tables \ref{tab: K 1e-6}-\ref{tab:a0=b0=c0=0}, are based on a setting where \(\nu=0.3\). These tables represent the situation where the thermo-poroelastic material is compressible. We also report the CPU times for the different algorithms. In particular, comparison between the coupled and decoupled algorithms reveals that the decoupled approaches incur lower computational costs, primarily because each subproblem in the decoupled scheme involves significantly smaller system sizes compared to the fully coupled formulation.

In Table \ref{tab: nu 0.49999}, we set a constant Poisson ratio of \(\nu = 0.49999\), while keeping all other physical parameters fixed. Since the Poisson ratio \(\nu\) is close to \(0.5\) which in other words means \( \lambda\to\infty\), the thermo-poroelastic material behaves nearly incompressibly, leading the mixed linear elasticity model to approach the incompressible Stokes model. The results presented in Table \ref{tab: nu 0.49999} demonstrate that the energy-norm errors for all algorithms exhibit optimal orders as the material approaches incompressibility, thereby emphasizing the robustness of these algorithms with respect to the parameter \(\nu\).

\subsection{Tests for the parameters $ K$ and $\Theta$}

In this subsection, we investigate the accuracy of the algorithms under varying settings of hydraulic conductivity \(\bm K\) and effective thermal conductivity \(\bm \Theta\). For reference, we have already tested the case \(\bm K = \bm \Theta = 0.1 \bm I\) in Table \ref{tab: K Theta 0.1 nu 0.3} in previous experiments. In the current tests (Table \ref{tab: K 1e-6}-Table \ref{tab:K Theta 1e-6}), we set \((\bm K, \bm\Theta) =(10^{-6}\bm I,0.1\bm I)\), \((\bm K, \bm\Theta) =(0.1\bm I, 10^{-6}\bm I)\), and \((\bm K, \bm\Theta) =(10^{-6}\bm I, 10^{-6}\bm I)\), respectively, to observe the impact of significantly smaller conductivity values on the solution accuracy.
The remaining key parameters are fixed as \(E = 1\), \(c_0 = a_0 = 0.2\), \(b_0 = 0.1\),  \(\alpha = \beta = 0.1\), and \(\nu=0.3\). We present the numerical results obtained using \textbf{Alg. 1} (the coupled algorithm) and \textbf{Alg. 2} (the iterative decoupling algorithm), with varying numbers of iterations, in Tables \ref{tab: K 1e-6} through \ref{tab:K Theta 1e-6}.

By comparing the results from Tables \ref{tab: K 1e-6} and \ref{tab:K Theta 1e-6} with those in Table \ref{tab: K Theta 0.1 nu 0.3}, we observe that all algorithms continue to produce solutions with optimal convergence rates. This indicates that the algorithms remain effective and robust, even when hydraulic and thermal conductivities are significantly reduced, further validating their performance under a wide range of physical conditions.

 \begin{table}[htbp]\tiny 
 \centering 
\caption{The errors and convergence orders with decreasing mesh sizes   $ \bm K=10^{-6}\bm{I}$. }
\begin{tabular}{ccccccccc}
\toprule
  $h$&$\Vert\bm u\Vert_{H^1}$ &rate &$\Vert\xi\Vert_{L^2}$ &rate &$\Vert p\Vert_{H^1}$ &rate&$\Vert T\Vert_{H^1}$ &rate \\ 
\midrule
Alg.1,\\
$\Delta t=10^{-3}$\\
\midrule
$1/16$ & 1.00615e-01 & 0.00 & 6.07391e-03 & 0.00 & 2.48018e-01 & 0.00 & 2.31172e-01 & 0.00 \\
$1/32$ & 2.53670e-02 & 1.99 & 1.50130e-03 & 2.02 & 1.13982e-01 & 1.12 & 1.09927e-01 & 1.07 \\
$1/64$ & 6.35858e-03 & 2.00 & 3.74103e-04 & 2.00 & 5.51891e-02 & 1.05 & 5.42244e-02 & 1.02 \\
$1/128$ & 1.59110e-03 & 2.00 & 9.32365e-05 & 2.00 & 2.72419e-02 & 1.02 & 2.70188e-02 & 1.00 \\
\midrule
Alg.2, \\
$\Delta t=5\times 10^{-3}$,\\
iters=5 &\\
\midrule
$1/16$ & 1.00615e-01 & 0.00 & 6.07747e-03 & 0.00 & 2.49194e-01 & 0.00 & 2.32137e-01 & 0.00 \\
$1/32$ & 2.53671e-02 & 1.99 & 1.50226e-03 & 2.02 & 1.14131e-01 & 1.13 & 1.10059e-01 & 1.08 \\
$1/64$ & 6.35859e-03 & 2.00 & 3.74397e-04 & 2.00 & 5.51695e-02 & 1.05 & 5.42417e-02 & 1.02 \\
$1/128$ & 1.59111e-03 & 2.00 & 9.33786e-05 & 2.00 & 2.72016e-02 & 1.02 & 2.70214e-02 & 1.01 \\
\midrule
Alg.2, \\
$\Delta t=10^{-2}$,\\
iters=10 &\\
\midrule
$1/16$ & 1.00616e-01 & 0.00 & 6.08155e-03 & 0.00 & 2.50703e-01 & 0.00 & 2.33492e-01 & 0.00 \\
$1/32$ & 2.53669e-02 & 1.99 & 1.50131e-03 & 2.02 & 1.14277e-01 & 1.13 & 1.10188e-01 & 1.08 \\
$1/64$ & 6.35831e-03 & 2.00 & 3.72183e-04 & 2.01 & 5.50831e-02 & 1.05 & 5.42265e-02 & 1.02 \\
$1/128$ & 1.59083e-03 & 2.00 & 9.11603e-05 & 2.03 & 2.71163e-02 & 1.02 & 2.70061e-02 & 1.01 \\
\bottomrule
\label{tab: K 1e-6} 
\end{tabular}  
\end{table}

\begin{table}[htbp]\tiny 
\centering
\caption{The errors and convergence orders with decreasing mesh sizes. $ \bm\Theta=10^{-6}\bm{I}$. }
\begin{tabular}{ccccccccc}
\toprule
  $h$&$\Vert\bm u\Vert_{H^1}$ &rate &$\Vert\xi\Vert_{L^2}$ &rate &$\Vert p\Vert_{H^1}$ &rate&$\Vert T\Vert_{H^1}$ &rate \\ 
\midrule
Alg.1, \\
$\Delta t=10^{-3}$\\
\midrule
$1/16$ & 1.00615e-01 & 0.00 & 6.07391e-03 & 0.00 & 2.31172e-01 & 0.00 & 2.48018e-01 & 0.00 \\
$1/32$ & 2.53670e-02 & 1.99 & 1.50130e-03 & 2.02 & 1.09927e-01 & 1.07 & 1.13982e-01 & 1.12 \\
$1/64$ & 6.35858e-03 & 2.00 & 3.74103e-04 & 2.00 & 5.42244e-02 & 1.02 & 5.51891e-02 & 1.05 \\
$1/128$ & 1.59110e-03 & 2.00 & 9.32365e-05 & 2.00 & 2.70188e-02 & 1.00 & 2.72419e-02 & 1.02 \\
\midrule
Alg.2,\\
$\Delta t=5\times 10^{-3}$,\\
iters=5 &\\
\midrule
$1/16$ & 1.00615e-01 & 0.00 & 6.07747e-03 & 0.00 & 2.32137e-01 & 0.00 & 2.49194e-01 & 0.00 \\
$1/32$ & 2.53671e-02 & 1.99 & 1.50226e-03 & 2.02 & 1.10059e-01 & 1.08 & 1.14131e-01 & 1.13 \\
$1/64$ & 6.35859e-03 & 2.00 & 3.74397e-04 & 2.00 & 5.42417e-02 & 1.02 & 5.51695e-02 & 1.05 \\
$1/128$ & 1.59111e-03 & 2.00 & 9.33786e-05 & 2.00 & 2.70214e-02 & 1.01 & 2.72016e-02 & 1.02 \\
\midrule
Alg.2,\\
$\Delta t=10^{-2}$,\\
iters=10 &\\
\midrule
$1/16$ & 1.00616e-01 & 0.00 & 6.08155e-03 & 0.00 & 2.33492e-01 & 0.00 & 2.50703e-01 & 0.00 \\
$1/32$ & 2.53669e-02 & 1.99 & 1.50131e-03 & 2.02 & 1.10188e-01 & 1.08 & 1.14277e-01 & 1.13 \\
$1/64$ & 6.35831e-03 & 2.00 & 3.72183e-04 & 2.01 & 5.42265e-02 & 1.02 & 5.50831e-02 & 1.05 \\
$1/128$ & 1.59083e-03 & 2.00 & 9.11603e-05 & 2.03 & 2.70061e-02 & 1.01 & 2.71163e-02 & 1.02 \\
\bottomrule
\label{tab:Theta 1e-6}   
\end{tabular}           
\end{table}

\begin{table}[htbp]  \tiny  
\centering 
\caption{The errors and convergence orders with decreasing mesh sizes. $ \bm K=10^{-6}\bm{I},\bm\Theta=10^{-6}\bm{I}$. }
\begin{tabular}{ccccccccc}
\toprule
  $h$&$\Vert\bm u\Vert_{H^1}$ &rate &$\Vert\xi\Vert_{L^2}$ &rate &$\Vert p\Vert_{H^1}$ &rate&$\Vert T\Vert_{H^1}$ &rate \\ 
\midrule
Alg.1, \\
$\Delta t=10^{-3}$\\
\midrule
$1/16$ & 1.00629e-01 & 0.00 & 6.18582e-03 & 0.00 & 2.73973e-01 & 0.00 & 2.73973e-01 & 0.00 \\
$1/32$ & 2.53705e-02 & 1.99 & 1.53143e-03 & 2.01 & 1.21482e-01 & 1.17 & 1.21482e-01 & 1.17 \\
$1/64$ & 6.35948e-03 & 2.00 & 3.81830e-04 & 2.00 & 5.71664e-02 & 1.09 & 5.71664e-02 & 1.09 \\
$1/128$ & 1.59132e-03 & 2.00 & 9.51228e-05 & 2.01 & 2.77233e-02 & 1.04 & 2.77233e-02 & 1.04 \\
\midrule
Alg.2, \\
$\Delta t=5\times 10^{-3}$, \\
iters=5 &\\
\midrule
$1/16$ & 1.00628e-01 & 0.00 & 6.18391e-03 & 0.00 & 2.73713e-01 & 0.00 & 2.73713e-01 & 0.00 \\
$1/32$ & 2.53705e-02 & 1.99 & 1.53140e-03 & 2.01 & 1.21507e-01 & 1.17 & 1.21507e-01 & 1.17 \\
$1/64$ & 6.35951e-03 & 2.00 & 3.82333e-04 & 2.00 & 5.72693e-02 & 1.09 & 5.72693e-02 & 1.09 \\
$1/128$ & 1.59137e-03 & 2.00 & 9.58070e-05 & 2.00 & 2.78597e-02 & 1.04 & 2.78597e-02 & 1.04 \\
\midrule
Alg.2,\\
$\Delta t=10^{-2}$,\\
iters=10 &\\
\midrule
$1/16$ & 1.00628e-01 & 0.00 & 6.18228e-03 & 0.00 & 2.73490e-01 & 0.00 & 2.73490e-01 & 0.00 \\
$1/32$ & 2.53701e-02 & 1.99 & 1.52778e-03 & 2.02 & 1.21011e-01 & 1.18 & 1.21011e-01 & 1.18 \\
$1/64$ & 6.35902e-03 & 2.00 & 3.78244e-04 & 2.01 & 5.67190e-02 & 1.09 & 5.67190e-02 & 1.09 \\
$1/128$ & 1.59093e-03 & 2.00 & 9.19929e-05 & 2.04 & 2.73453e-02 & 1.05 & 2.73453e-02 & 1.05 \\
\bottomrule
\label{tab:K Theta 1e-6}
\end{tabular}    
\end{table}

\subsection{ Tests for the parameter $a_0,b_0,c_0$}

\begin{table}[htbp]\tiny 
\centering 
\caption{ The errors and convergence orders with decreasing mesh sizes. $a_0=b_0= c_0=0$. } 
\begin{tabular}{ccccccccc}
\toprule
  $h$&$\Vert\bm u\Vert_{H^1}$ &orders &$\Vert\xi\Vert_{L^2}$ &orders &$\Vert p\Vert_{H^1}$ &orders&$\Vert T\Vert_{H^1}$ &orders \\ 
\midrule
Alg.1,\\
$\Delta t=10^{-3}$\\
\midrule
$1/16$ & 1.00716e-01 & 0.00 & 6.74538e-03 & 0.00 & 2.60803e-01 & 0.00 & 2.60803e-01 & 0.00 \\
$1/32$ & 2.53932e-02 & 1.99 & 1.67613e-03 & 2.01 & 1.14106e-01 & 1.19 & 1.14106e-01 & 1.19 \\
$1/64$ & 6.36502e-03 & 2.00 & 4.17182e-04 & 2.01 & 5.47500e-02 & 1.06 & 5.47500e-02 & 1.06 \\
$1/128$ & 1.59255e-03 & 2.00 & 1.02942e-04 & 2.02 & 2.70774e-02 & 1.02 & 2.70774e-02 & 1.02 \\
\midrule
Alg.2, \\
$\Delta t=5\times 10^{-3}$,\\
iters=5 &\\
\midrule
$1/16$ & 1.00807e-01 & 0.00 & 7.31400e-03 & 0.00 & 2.91681e-01 & 0.00 & 2.91681e-01 & 0.00 \\
$1/32$ & 2.54374e-02 & 1.99 & 1.95023e-03 & 1.91 & 1.22342e-01 & 1.25 & 1.22342e-01 & 1.25 \\
$1/64$ & 6.40628e-03 & 1.99 & 6.41686e-04 & 1.60 & 5.85438e-02 & 1.06 & 5.85438e-02 & 1.06 \\
$1/128$ & 1.66723e-03 & 1.94 & 3.49812e-04 & 0.88 & 3.04263e-02 & 0.94 & 3.04263e-02 & 0.94 \\
\midrule
Alg.2, \\
$\Delta t=10^{-2}$,\\
iters=10 &\\
\midrule
$1/16$ & 1.00916e-01 & 0.00 & 7.96738e-03 & 0.00 & 3.30774e-01 & 0.00 & 3.30774e-01 & 0.00 \\
$1/32$ & 2.54412e-02 & 1.99 & 1.97495e-03 & 2.01 & 1.24587e-01 & 1.41 & 1.24587e-01 & 1.41 \\
$1/64$ & 6.37477e-03 & 2.00 & 4.79389e-04 & 2.04 & 5.59219e-02 & 1.16 & 5.59219e-02 & 1.16 \\
$1/128$ & 1.59309e-03 & 2.00 & 1.06986e-04 & 2.16 & 2.71232e-02 & 1.04 & 2.71232e-02 & 1.04 \\
\bottomrule
 \label{tab:a0=b0=c0=0}
\end{tabular}           
\end{table} 
In this subsection, we investigate the influence of the effective thermal parameter \(a_0\), the thermal dilation coefficient \(b_0\), and the specific storage coefficient \(c_0\) on the accuracy of the algorithms. As discussed in Remark \ref{re: a0-b0=c0-b0}, the case \(c_0 - b_0 = a_0 - b_0 = 0\) may affect the convergence rate of the iterative decoupling algorithms. To evaluate this, Table \ref{tab:a0=b0=c0=0} presents the limiting scenario where the coefficients \(a_0\), \(b_0\), and \(c_0\) are all set to zero. For the remaining parameters, we fix \(E = 1\), \(\alpha = \beta = 0.1\), and \(\bm{K} = \bm{\Theta} = 0.1\bm{I}\). 

A comparison of the second entry in Table \ref{tab:a0=b0=c0=0} with the corresponding entry in Table \ref{tab: K Theta 0.1 nu 0.3} reveals that setting \(a_0 = b_0 = c_0 = 0\) leads to a slight degradation in the errors and convergence rates for the iterative decoupling algorithm (\textbf{Alg. 2}). However, when the number of iterations increases to \(10\), errors and rates reach optimal levels, as shown in the third entry of Table \ref{tab:a0=b0=c0=0}. These results confirm that the energy norm errors for all algorithms maintain optimal convergence orders, underscoring their robustness with respect to the parameters \(a_0\), \(b_0\), and \(c_0\).

\subsection{Tests with a large terminal time $\tau$}
In this subsection, we evaluate the performance of the proposed algorithms for a relatively long terminal time $\tau = 1$. Owing to the extended simulation duration, all algorithms require a greater number of time steps to complete the computation. The parameters are set as follows: $\nu = 0.3$, $E = 1$, $a_0 = c_0 = 0.2$, $b_0 = 0.1$, $\alpha = \beta = 0.1$, and $\bm{K} = \bm{\Theta} = 0.1 I$.

Table~\ref{tab: same Dt of Alg1 and Alg2} presents a comparison of \textbf{Alg. 1} and \textbf{Alg. 2} under different time-step sizes, specifically $\Delta t = 10^{-3}$, $\Delta t = 5 \times 10^{-3}$, and $\Delta t = 10^{-2}$. 
In all these cases, the spatial discretization error remains dominant. The numerical results show that both the coupled and decoupled algorithms yield comparable accuracy. As expected, reducing the time-step size leads to improved accuracy for both algorithms. Furthermore, for fixed time-step sizes, increasing the number of iterations in \textbf{Alg. 2} from 5 to 10 results in errors that more closely align with those of \textbf{Alg. 1}, thereby confirming the theoretical convergence behavior stated in Theorem~\ref{th: convergence_of_N1}.

As in the comparisons presented in Table~\ref{tab: K Theta 0.1 nu 0.3}, we examine \textbf{Alg. 1} using a time-step size of $\Delta t = 1.0 \times 10^{-3}$, and \textbf{Alg. 2} using $\Delta t = 5 \times 10^{-3}$ with 5 iterations or $\Delta t = 10^{-2}$ with 10 iterations. The results show that the numerical errors for both algorithms are comparable. These observations confirm that the iterative decoupled algorithm remains accurate and robust even when applied to simulations with large terminal times. Moreover, consistent with the findings in Table~\ref{tab: K Theta 0.1 nu 0.3}, the decoupled algorithm offers improved computational efficiency.

\begin{table}[htbp] \tiny  
\caption{The errors and convergence orders of different algorithms. $\tau=1$.}   
\centering
\begin{tabular}{ccccccccc}
\toprule
  $h$&$\Vert\bm u\Vert_{H^1}$ &orders &$\Vert\xi\Vert_{L^2}$ &orders &$\Vert p\Vert_{H^1}$ &orders&$\Vert T\Vert_{H^1}$ &orders \\ 
\midrule
Alg.1,$\Delta t=10^{-2}$&\\
\midrule
$1/8$ & 1.45210e-01 & 0.00 & 9.12749e-03 & 0.00 & 1.58002e-01 & 0.00 & 1.58002e-01 & 0.00 \\
$1/16$ & 3.73731e-02 & 1.96 & 2.16086e-03 & 2.08 & 7.99221e-02 & 0.98 & 7.99221e-02 & 0.98 \\
$1/32$ & 9.42251e-03 & 1.99 & 5.33182e-04 & 2.02 & 4.00855e-02 & 1.00 & 4.00855e-02 & 1.00 \\
$1/64$ & 2.36208e-03 & 2.00 & 1.34243e-04 & 1.99 & 2.00727e-02 & 1.00 & 2.00727e-02 & 1.00 \\
\midrule
Alg.2,$\Delta t= 10^{-2}$,iters=5 &\\
\midrule
$1/8$ & 1.45210e-01 & 0.00 & 9.12737e-03 & 0.00 & 1.58000e-01 & 0.00 & 1.58000e-01 & 0.00 \\
$1/16$ & 3.73731e-02 & 1.96 & 2.16072e-03 & 2.08 & 7.99198e-02 & 0.98 & 7.99198e-02 & 0.98 \\
$1/32$ & 9.42248e-03 & 1.99 & 5.32938e-04 & 2.02 & 4.00813e-02 & 1.00 & 4.00813e-02 & 1.00 \\
$1/64$ & 2.36200e-03 & 2.00 & 1.33584e-04 & 2.00 & 2.00645e-02 & 1.00 & 2.00645e-02 & 1.00 \\
\midrule
Alg.2,$\Delta t= 10^{-2}$,iters=10 &\\
\midrule
$1/8$ & 1.45210e-01 & 0.00 & 9.12749e-03 & 0.00 & 1.58002e-01 & 0.00 & 1.58002e-01 & 0.00 \\
$1/16$ & 3.73731e-02 & 1.96 & 2.16086e-03 & 2.08 & 7.99221e-02 & 0.98 & 7.99221e-02 & 0.98 \\
$1/32$ & 9.42251e-03 & 1.99 & 5.33182e-04 & 2.02 & 4.00855e-02 & 1.00 & 4.00855e-02 & 1.00 \\
$1/64$ & 2.36208e-03 & 2.00 & 1.34243e-04 & 1.99 & 2.00727e-02 & 1.00 & 2.00727e-02 & 1.00 \\
\midrule
\midrule
Alg.1,$\Delta t=5\times 10^{-3}$&\\
\midrule
$1/8$ & 1.45210e-01 & 0.00 & 9.12771e-03 & 0.00 & 1.58002e-01 & 0.00 & 1.58002e-01 & 0.00 \\
$1/16$ & 3.73731e-02 & 1.96 & 2.16106e-03 & 2.08 & 7.99192e-02 & 0.98 & 7.99192e-02 & 0.98 \\
$1/32$ & 9.42250e-03 & 1.99 & 5.33110e-04 & 2.02 & 4.00783e-02 & 1.00 & 4.00783e-02 & 1.00 \\
$1/64$ & 2.36193e-03 & 2.00 & 1.33096e-04 & 2.00 & 2.00576e-02 & 1.00 & 2.00576e-02 & 1.00 \\
\midrule
Alg.2,$\Delta t=5\times 10^{-3}$,iters=5 &\\
\midrule
$1/8$ & 1.45210e-01 & 0.00 & 9.12758e-03 & 0.00 & 1.58001e-01 & 0.00 & 1.58001e-01 & 0.00 \\
$1/16$ & 3.73731e-02 & 1.96 & 2.16092e-03 & 2.08 & 7.99180e-02 & 0.98 & 7.99180e-02 & 0.98 \\
$1/32$ & 9.42248e-03 & 1.99 & 5.32940e-04 & 2.02 & 4.00764e-02 & 1.00 & 4.00764e-02 & 1.00 \\
$1/64$ & 2.36189e-03 & 2.00 & 1.32775e-04 & 2.00 & 2.00541e-02 & 1.00 & 2.00541e-02 & 1.00 \\
\midrule
Alg.2,$\Delta t= 5 \times 10^{-3}$,iters=10 &\\
\midrule
$1/8$ & 1.45210e-01 & 0.00 & 9.12771e-03 & 0.00 & 1.58002e-01 & 0.00 & 1.58002e-01 & 0.00 \\
$1/16$ & 3.73731e-02 & 1.96 & 2.16106e-03 & 2.08 & 7.99192e-02 & 0.98 & 7.99192e-02 & 0.98 \\
$1/32$ & 9.42250e-03 & 1.99 & 5.33110e-04 & 2.02 & 4.00783e-02 & 1.00 & 4.00783e-02 & 1.00 \\
$1/64$ & 2.36193e-03 & 2.00 & 1.33096e-04 & 2.00 & 2.00576e-02 & 1.00 & 2.00576e-02 & 1.00 \\
\midrule
\midrule
Alg.1, $\Delta t=10^{-3}$ &\\
\midrule
$1/8$ & 1.45210e-01 & 0.00 & 9.12790e-03 & 0.00 & 1.58003e-01 & 0.00 & 1.58003e-01 & 0.00 \\
$1/16$ & 3.73731e-02 & 1.96 & 2.16126e-03 & 2.08 & 7.99186e-02 & 0.98 & 7.99186e-02 & 0.98 \\
$1/32$ & 9.42251e-03 & 1.99 & 5.33223e-04 & 2.02 & 4.00762e-02 & 1.00 & 4.00762e-02 & 1.00 \\
$1/64$ & 2.36189e-03 & 2.00 & 1.32863e-04 & 2.00 & 2.00529e-02 & 1.00 & 2.00529e-02 & 1.00 \\
\midrule
Alg.2,$\Delta t= 10^{-3}$,iters=5 &\\
\midrule
$1/8$ & 1.45210e-01 & 0.00 & 9.12776e-03 & 0.00 & 1.58002e-01 & 0.00 & 1.58002e-01 & 0.00 \\
$1/16$ & 3.73731e-02 & 1.96 & 2.16113e-03 & 2.08 & 7.99189e-02 & 0.98 & 7.99189e-02 & 0.98 \\
$1/32$ & 9.42250e-03 & 1.99 & 5.33150e-04 & 2.02 & 4.00774e-02 & 1.00 & 4.00774e-02 & 1.00 \\
$1/64$ & 2.36191e-03 & 2.00 & 1.32987e-04 & 2.00 & 2.00556e-02 & 1.00 & 2.00556e-02 & 1.00 \\
\midrule
Alg.2,$\Delta t=10^{-3}$,iters=10 &\\
\midrule
$1/8$ & 1.45210e-01 & 0.00 & 9.12790e-03 & 0.00 & 1.58003e-01 & 0.00 & 1.58003e-01 & 0.00 \\
$1/16$ & 3.73731e-02 & 1.96 & 2.16126e-03 & 2.08 & 7.99186e-02 & 0.98 & 7.99186e-02 & 0.98 \\
$1/32$ & 9.42251e-03 & 1.99 & 5.33223e-04 & 2.02 & 4.00762e-02 & 1.00 & 4.00762e-02 & 1.00 \\
$1/64$ & 2.36189e-03 & 2.00 & 1.32863e-04 & 2.00 & 2.00529e-02 & 1.00 & 2.00529e-02 & 1.00 \\
\bottomrule
\label{tab: same Dt of Alg1 and Alg2}
\end{tabular}       
\end{table}

\section{Conclusions}
In this paper, we propose a novel algorithm designed to decouple the computations of the linear thermo-poroelastic model, based on a newly developed four-field formulation. The algorithm effectively separates the numerical computations into two sub-models, streamlining the overall computational process. We perform a comprehensive error analysis to establish the algorithm's unconditional stability and optimal convergence. Furthermore, numerical experiments are performed to validate the theoretical error estimates. The results demonstrate that the proposed algorithm is unconditionally stable, computationally efficient, and achieves optimal convergence rates even under widely varying parameter settings or in the limiting cases of parameters, highlighting its robustness and reliability for addressing thermo-poroelastic problems.

\section*{Acknowledgments}
The authors thank two reviewers for their constructive suggestions and comments. 
The work of M. Cai is partially supported by Army Research Office award W911NF-23-1-0004, the NIH-RCMI award (Grant No. 347U54MD013376), and the affiliated project award from the Center for Equitable Artificial Intelligence and Machine Learning Systems (CEAMLS) at Morgan State University (project ID 02232301). The work of J. Li and Z. Li are partially supported by the Shenzhen Sci-Tech Fund No. RCJC20200714114556020, 
Guangdong Basic and Applied Research Fund No. 2023B1515250005.
The work of Q. Liu is partially supported by the NSFC (12271186, 12271178), the Guangdong Basic and Applied Basic Research Foundation (2022B1515120009), and the Science and Technology Program of Shenzhen, China(20231121110406001).

\bibliographystyle{siamplain}
\bibliography{ref}

\end{document}